\documentclass{article}
\usepackage[T1]{fontenc}
\usepackage[utf8]{inputenc}
\usepackage{multirow} 
\usepackage{amsfonts,mathtools,amssymb}
\usepackage{graphicx,wrapfig,tikz,animate,subcaption}
\usepackage{lmodern}
\usepackage{fancyhdr,graphicx,wrapfig,amsmath,amssymb, mathtools, scrextend, titlesec, enumitem,color, soul,xcolor}
\usepackage{array}
\usepackage{hyperref}
\usepackage{caption}
\usepackage{comment}
\usepackage{cancel}
\usepackage{tabularx}
\DeclareMathOperator{\Pe}{Pe}

\usepackage[boxed,linesnumbered]{algorithm2e} 

\usetikzlibrary{arrows}

\definecolor{verylightgray}{gray}{0.93}
\usepackage[color=verylightgray]{todonotes}

\usepackage{wrapfig,tikz}
\usetikzlibrary{shapes.misc}

\newcommand{\Change}[2]{{\delete{#1 }}\new{#2}}
\newcommand{\delete}[1]{}
\newcommand{\new}[1]{{#1}}
\newtheorem{thm}{Teorema}

\newtheorem{rem}[thm]{Remark}
\usepackage{xspace}
\newcommand{\N}{\ensuremath{\mathbb{N}}\xspace}
\newcommand{\R}{\ensuremath{\mathbb{R}}\xspace}
\newcommand{\PP}{\ensuremath{\mathcal{P}}\xspace}
\newcommand{\Q}{\ensuremath{\mathcal{Q}}\xspace}
\newcommand{\T}{\ensuremath{\mathcal{T}}\xspace}
\newcommand{\aaa}{\ensuremath{\mathbf{a}}\xspace}
\newcommand{\nnn}{\ensuremath{\mathbf{n}}\xspace}
\newcommand{\Vrfbn}{\ensuremath{V_{B\!M\!Z}}\xspace}

\newcommand{\inflow}{\partial_{\textbf{-}}}
\DeclareMathOperator{\Div}{div}
\renewcommand{\span}{\text{span}}

\newcommand{\qtfa}{\quad\text{for all }}
\newcommand{\tand}{\quad\text{and}\quad}
\newcommand{\Hoi}[1][\Omega]{\ensuremath{H_0^1(#1)}\xspace}
\newcommand{\Hgi}[1][\Omega]{\ensuremath{H_g^1(#1)}\xspace}
\newcommand{\bT}[1][{}]{\ensuremath{\psi_T^{#1}}}
\newcommand{\tbT}[1][{}]{\ensuremath{\tilde \psi_T^{#1}}}
\newcommand{\bS}[1][{}]{b_S^{#1}}
\newcommand{\patchS}[1][{}]{\omega_S^{#1}}

\title{Patch bubbles for advection-dominated steady and unsteady problems}
\author{Eberhard B\"ansch\footnote{Applied Mathematics III,
Friedrich-Alexander-University Erlangen-Nuremberg,
Cauerstr. 11, 91058 Erlangen, Germany
\href{mailto:baensch@math.fau.de}{baensch@math.fau.de}}
\and 
Pedro Morin\footnote{Universidad Nacional del Litoral and CONICET, 
Departamento de Matem\'atica, 
Facultad de Ingenier\'{\i}a Qu\'{\i}mica, Santiago del Estero 2829, S3000AOM Santa Fe, Argentina. 
\href{mailto:pmorin@fiq.unl.edu.ar}{pmorin@fiq.unl.edu.ar}}
\and Itatí Zocola\footnote{Universidad Nacional del Litoral and CONICET, 
Departamento de Matem\'atica, 
Facultad de Ingenier\'{\i}a Qu\'{\i}mica, Santiago del Estero 2829, S3000AOM Santa Fe, Argentina. 
\href{mailto:izocola@fiq.unl.edu.ar}{izocola@fiq.unl.edu.ar}}
}
\date{}

\begin{document}
\hyphenpenalty = 5000

	\maketitle


\begin{abstract}
A novel variant of the \emph{residual-free bubble} method (RFB) for advection-dominated problems is presented. Since the usual RFB still suffers from oscillations and strong under/overshoots, the bubble space is enriched by \emph{patch bubbles}, giving more freedom to the bubble space.

We use a recursive and efficient approach to accurately compute the bubbles. Numerical experiments clearly demonstrate the superiority of our method compared to the standard RFB.

While similar in spirit to the \emph{enhanced residual-free bubble} (eRFB) method by Cangiani and S\"uli~\cite{Cangiani2005long,Cangiani2005short}, our approach differs in the definition and computation of the additional bubbles.

In addition, we extend the methodology to solve problems with non constant coefficients and we develop a novel bubble-based stabilization technique for time-dependent problems, which performs very accurately.
\end{abstract}


\tableofcontents

\section{Introduction}

In this article we consider the following advection-reaction-diffusion problem, in a domain $\Omega$ of $\R^2$:
\begin{equation*}
\left\{
\begin{aligned}
    u_t - \epsilon \Delta u + \aaa\cdot \nabla  u + c \, u &=f ,\qquad& &\text{in $\Omega$, $t>0$,}\\
        u &= g, \qquad& &\text{on $\partial\Omega$, } t>0,\\
    u(\cdot,0) &= u^0, \qquad& & \text{in }\Omega.\\
\end{aligned}
    \right.
\end{equation*}
Here
$\epsilon>0$,
$\aaa \in W^{1,\infty}(\Omega)$, $c\in L^\infty(\Omega)$, 
with 
$c-\frac12 \Div \aaa \ge 0$,
$u^0 \in L^2(\Omega)$, and
$f(t)\in L^2(\Omega)$, 
$g(t) \in H^{1/2}(\partial\Omega)$, for all $t>0$.

This problem is called \emph{advection dominated} when $\epsilon \ll |\aaa|$, \Change{}{or equivalently, the Péclet number $\Pe:=|\aaa|/\epsilon\gg1$}.
It is well known that in this case the standard Galerkin method using polynomial finite elements for both the trial and test functions suffer from instabilities that manifest in high oscillations on the discrete solution, in the steady as well as the unsteady case.
In order to cure this phenomenon and make finite element methods applicable to advection dominated problems, several attempts have been made.

    The literature extensively documents methods for solving advection-diffusion-reaction (ADR) equations using stabilized finite element techniques. 
    A detailed analysis in {\cite{Codina1998}} focused on various stabilized approaches, clarifying similarities and differences, including those of the Streamline Upwind Petrov-Galerkin (SUPG) method. Parameters suitable for SUPG on steady-state ADR equations were investigated in {\cite{Codina1998}}. Additionally, specialized stabilization methods, like the Unusual Stabilized Finite Element Method (USFEM) {\cite{FrancaValentin2000}} have been explored. 
    Attempts to manage unwanted effects, such as a non-linear variant of SUPG aiming to suppress spurious oscillations 
    were studied in {\cite{Knopp2002}} and later extended to higher-order elements {\cite{lube2006residual}}. Other work has focused on methods for time-dependent ADR equations rooted in variational multiscale principles {\cite{CodinaBlasco2002},\cite{HaukeDoweidar2005},\cite{VolkerEtAl2006},\cite{GravemeierWall2008}}. 
    It has long been recognized that the SUPG method can yield oscillatory solutions, even for simpler steady-state convection-diffusion equations. 
    To counter this, a suite of techniques known as Spurious Oscillations at Layers Diminishing (SOLD) methods were introduced {\cite{john2007comparison},\cite{John2008}}. 
    Despite their success in steady-state contexts, very few $\text{SOLD}$ schemes have been successfully applied to the more complex time-dependent simulations {\cite{Bazilevs_2007_YZbeta}}. 
        \new{A change of variable approach and discontinuity capturing methods were proposed in~{\cite{real-world-problems}} to ensure physical constraints for advection-diffusion-reaction equations.
    A two-level finite element method for convection-diffusion-reaction problems, which consists of a global problem and an ensemble of local problems has been studied in~\cite{kaya-braack}.}
Most of the aforementioned stabilization methods enhance the standard Galerkin finite element formulation by explicitly adding extra terms to the differential equation. 
    Another alternative that has been explored is the use of Discontinuous Galerkin methods (DG) {\cite{Cockburn1999}}; these methods are particularly attractive since they achieve formal high-order accuracy and nonlinear stability while maintaining the ability to capture the strong gradients without producing strong spurious oscillations.

Another possibility to cure these instabilities is the use of \emph{residual-free bubbles}, which were first introduced for the steady case in~\cite{Brezzi1994}, and later analyzed in~\cite{Brezzi1999}, where the so-called \emph{stability norm} was introduced for measuring the error.
Among other things it was proved that this method coincides with the famous SUPG method~\cite{BFHR1997, BrooksHughes:82}, providing with a precise proposal for the parameters involved in the latter \new{(see Remark~\ref{rem:SUPG=RFB} below)}. 
Most of the other stabilization methods described in the previous paragraph require the tuning of \Change{unknown}{some parameters, when simulating real-world problems with a broad range of flow conditions, such as those considered in~\cite{real-world-problems}. For the simple geometries considered in the examples presented in this article there are some canonical choices which work very well in practice (see~\cite[Section 3.2]{JohnSchmeyer}). 
However, residual-free bubbles are attractive because they do not seem to need such a tuning, and that is why we have chosen to investigate further this kind of methods.}


\delete{
When trying to solve unsteady problems using residual-free bubbles for stabilization, we were not satisfied by the observed numerical results and thus embarked into a project trying to improve the method by adding bubbles supported not only on elements, but also on patches of two elements. 
In the meantime we came across the article~\cite{Cangiani2005long}, where bubbles supported on patches of two elements were first proposed, and its approximation properties were analyzed, again for steady problems. Their idea differs from ours, and we will explain the details later, but it is worth mentioning here that numerically their method and ours perform very similarly.
}
A thorough overview of residual-free bubble methods for steady advection-dominated problems can be found in~\cite{Cangiani2005long}, we thus restrict the list of references to those strictly related to our work.

\new{The main goal of this article is to propose and investigate the performance of a new residual free bubble method for advection-dominated advection-reaction-diffusion problems, which does not exhibit spurious oscillations, and is suitable for steady and also for unsteady problems. }
We do this in three steps.
First, we incorporate residual-free bubbles supported on patches of two elements, in order to add more freedom to the method across interelement edges; we do this for the steady case. 
Secondly, since the residual-free bubbles are solutions to a problem as complicated as the original one, we propose to use a recursive algorithm to compute the bubbles. With this recursion, the algorithm needs to compute $O(\hspace{0.05 cm} \log (\Change{\hspace{0.05 cm} |\aaa| h \hspace{0.05 cm}/\hspace{0.05 cm}\epsilon}{\Pe_h} ))$ bubbles, \Change{}{where $\Pe_h:=|\aaa| h /\epsilon$ is the element Péclet number}; the one corresponding to the deepest recursion level does not need stabilization. 
Lastly, we extend the method to the unsteady case.
\new{Recently, four space-time finite element methods have been proposed {\cite{space-time-AD}} for the unsteady advection-diffusion equation. In this article, instead, we focus on a time-stepping method, which is the method of choice for many practitioners, see~\cite{Burman} and the references therein.}

The outline of this article is as follows.
In Section~\ref{S:RFB-stat} we describe the original residual-free bubble method for steady problems, with bubbles supported on elements.
In Section~\ref{S:RFB-patch-stat} we describe our proposal to use residual-free bubbles supported on patches of two elements, we describe the recursive algorithm and implementation details. We end this section with some numerical experiments and comparisons with other methods.
In Section~\ref{S:RFB-instat} we present an algorithm for solving unsteady advection-dominated advection reaction diffusion problems, using these ideas together with some numerical experiments, to illustrate its performance.
We end this article with some conclusions and remarks about future lines of research.

\section{Residual-Free Element Bubbles for steady Problems}
\label{S:RFB-stat}

The first ideas about residual-free bubbles were presented for the steady problem:
\begin{equation}\label{eq:convdiff-stat-strong}
    \left\{
\begin{aligned}
    - \epsilon \Delta u + \aaa\cdot \nabla  u  + c\, u &=f \qquad& &\text{in $\Omega$,}\\
        u &= g, \qquad& &\text{on $\partial\Omega$} 
\end{aligned}
    \right.
\end{equation}
with (piecewise) constant $\aaa$, and $f$, vanishing reaction term $c=0$, and homogeneous boundary condition $g=0$. 

The associated bilinear form becomes:
\begin{equation}\label{eq:bilinear-op-a}
    a(\psi,\varphi) = \int_\Omega \epsilon \nabla  \psi \cdot \nabla  \varphi + \aaa \cdot \nabla  \psi \, \varphi + c \, \psi \, \varphi
\end{equation}
and the weak formulation of~\eqref{eq:convdiff-stat-strong} reads:
\[
\text{Find } u \in \Hoi : \qquad a(u,\varphi) = (f,\varphi) \qtfa \varphi \in \Hoi,
\]
where $(f,\varphi)=\int_\Omega f\varphi$, as usual, 
and $\Hoi = \{ v \in H^1(\Omega) : v_{|\partial \Omega} = 0\}$.

\Change{}{
A standard criterion to quantify the relative strength of advection and diffusion is given by the Péclet number. Considering a partition with diameter $h$, the element Péclet number is defined as $\Pe_h:=|a|h/\epsilon$. The problem, relative to the partition, is advection-dominated when $\Pe_h \gg 1$,  while it becomes diffusion-dominated for $\Pe_h \ll 1$. In the first case, sharp boundary layers may appear at a scale finer than the mesh-size and standard finite element discretizations typically exhibit spurious oscillations.}


In this section we recall the definition of the residual-free bubbles initially proposed in~\cite{Brezzi1994} for simplicial meshes with linear finite elements, assuming the reaction term to be zero.

Consider a bounded domain $\Omega \subset \mathbb{R}^d$ ($d=2,3$) and a conforming triangulation of simplices ${\T_h}$ of $\bar{\Omega}$, 
the finite element space $V_L$ is defined as
\[
V_L = \{ v \in \Hoi : v_{|T} \in \PP^1, \forall T\in \T_h\},
\]
with $\PP^1$ the space of polynomials of degree less than or equal to $1$.

For a triangulation of simplices $\T_h$ and for each element $T\in \T_h$ they defined the \emph{residual-free bubble} $\bT$ (denoted by $b_1^T$ in~\cite{Brezzi1994, Brezzi1999}) as follows:
\[
\bT\in \Hoi[T]:\qquad a(\bT,\varphi)=(1,\varphi) \qtfa \varphi\in \Hoi[T].
\]

Then, for every $T\in\T_h$, consider the one-dimensional space $B_T$ spanned by $\bT$  and set 
\[
V_B = \bigoplus_{T\in\T}B_T
\tand
V_h = V_L \oplus V_B.
\]
The discrete problem now reads
\begin{equation}\label{eq:DiscreteProblemBrezzi}
\text{Find } u_h = u_L + u_B \in V_h:\qquad
a(u_h,v_h) = (f,v_h)\qtfa v_h \in V_h. 
\end{equation}

Error estimates for the residual-free bubble method are provided in \cite[Theorem 4.2]{Brezzi1999}. Nevertheless, these bounds depend on the $H^s$ regularity of the solution (with $1 < s \leq 2$), which typically behaves as $1/\varepsilon^s$. As a consequence, the resulting estimates are not entirely satisfactory for convection-dominated regimes. The main result is established in the so-called stability norm, defined as:
\begin{equation}\label{def:stab-norm}
    \left(\epsilon ||\nabla  v||_{L^2(\Omega)}^2 + \sum_T h_T||\aaa \cdot \nabla  v||_{L^2(T)}^2\right)^{1/2}.
\end{equation}

\begin{rem}[Static condensation and SUPG]\label{rem:SUPG=RFB}\rm
After testing separately for $v_h \in V_B$ and $v_h \in V_L$ one can perform a static condensation 
and arrive at a linear system of equations for $u_L$, which is what practitioners are usually interested in.
    In this respect, the use of bubbles can be regarded as a \emph{stabilization method} for the finite element solution $u_L$, and coincides with the famous SUPG method~, with a precise definition of the stabilization constants. 
    \new{More precisely, in~{\cite[eq.~(43)]{BFHR1997}} the SUPG method is stated with a stabilization parameter $\tau_T = \frac{h_T}{2|\aaa|}$ for each element $T$ with $\Pe_T\ge 1$, and it is shown that the RFB method coincides with SUPG with stabilization parameter $\hat \tau_T = \frac{1}{|T|}\int_T b_T$ with $b_T$ the residual free bubble on $T$, i.e., $-\epsilon\Delta b_T + \aaa\cdot\nabla b_T = 1$ in $T$, $b_T = 0$ on $\partial T$.}   
\end{rem}

\begin{rem}[Computation of terms involving $\bT$]\rm
    It is worth noticing that integration by parts leads to
    \begin{equation}\label{eq:bT-identity-1}
    \epsilon \| \nabla  \bT \|_{0,T}^2  = \int_T \bT,
    \end{equation}
    and, for every affine function $\varphi$ on $T$,
    \begin{equation}\label{eq:bT-identity-2}
        \begin{aligned}
            \int_T \epsilon \nabla \bT \cdot\nabla  \varphi &= 0 \tand \\
            \int_T \aaa \cdot \nabla  \bT \, \varphi 
            &= - \int_T \aaa \cdot \nabla  \varphi \, \bT
            = -\aaa\cdot\nabla  \varphi_{|T}\int_T \bT.
        \end{aligned}
        \end{equation}
    Hence all the terms in the linear system involving $\bT$ can be approximated by its integral, and no other features of $\bT$. In the advection dominated regime, \Change{}{$\Pe_h \gg 1$}, this integral can be well approximated by the integral of the \emph{pyramid function} $\tbT$ which solves the pure advection problem
    \[
    \aaa\cdot \nabla  \tbT = 1, \ \text{in $T$}, 
    \qquad
    \tbT = 0, \ \text{on $\inflow T$},
    \]
    and can be computed exactly; here $\inflow T$ denotes the \emph{inflow} part of $\partial T$ (where $\aaa \cdot \nnn\le 0$).
\end{rem}

\begin{figure}[h!tbp]
    
    \hfil
    \includegraphics[width=.49\textwidth]{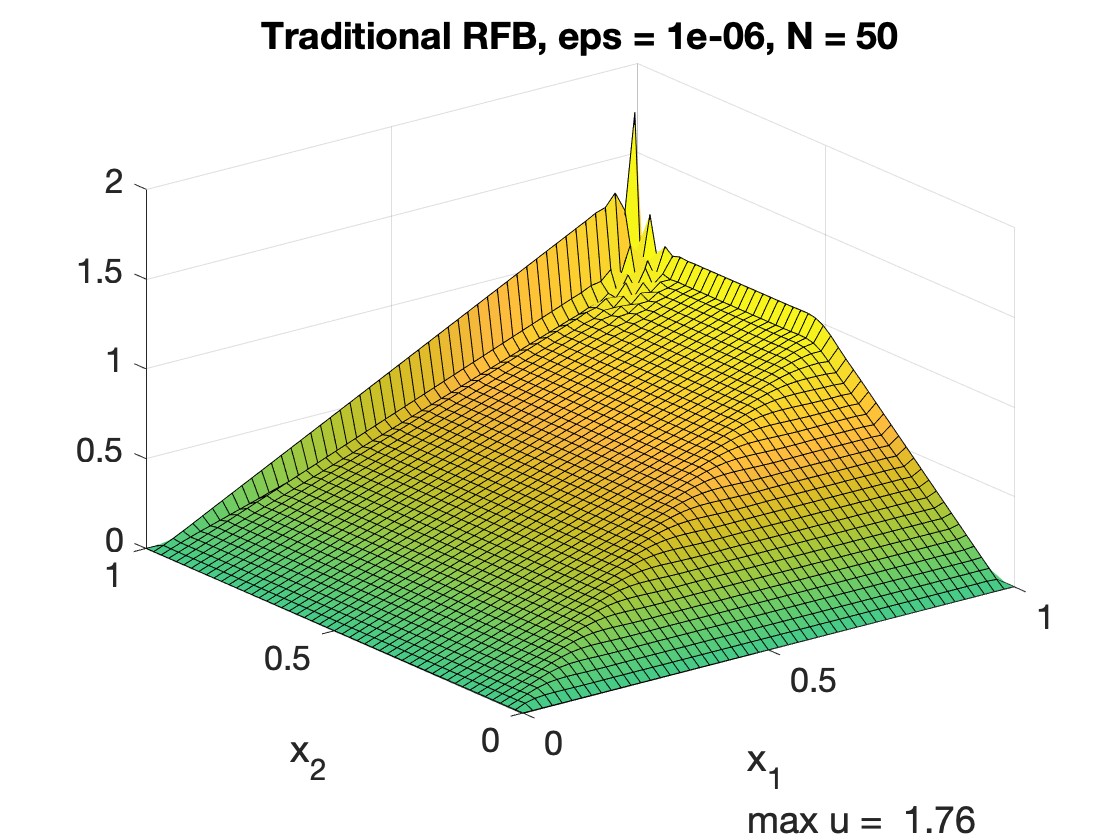}
    \hfil
    \includegraphics[width=.49\textwidth]{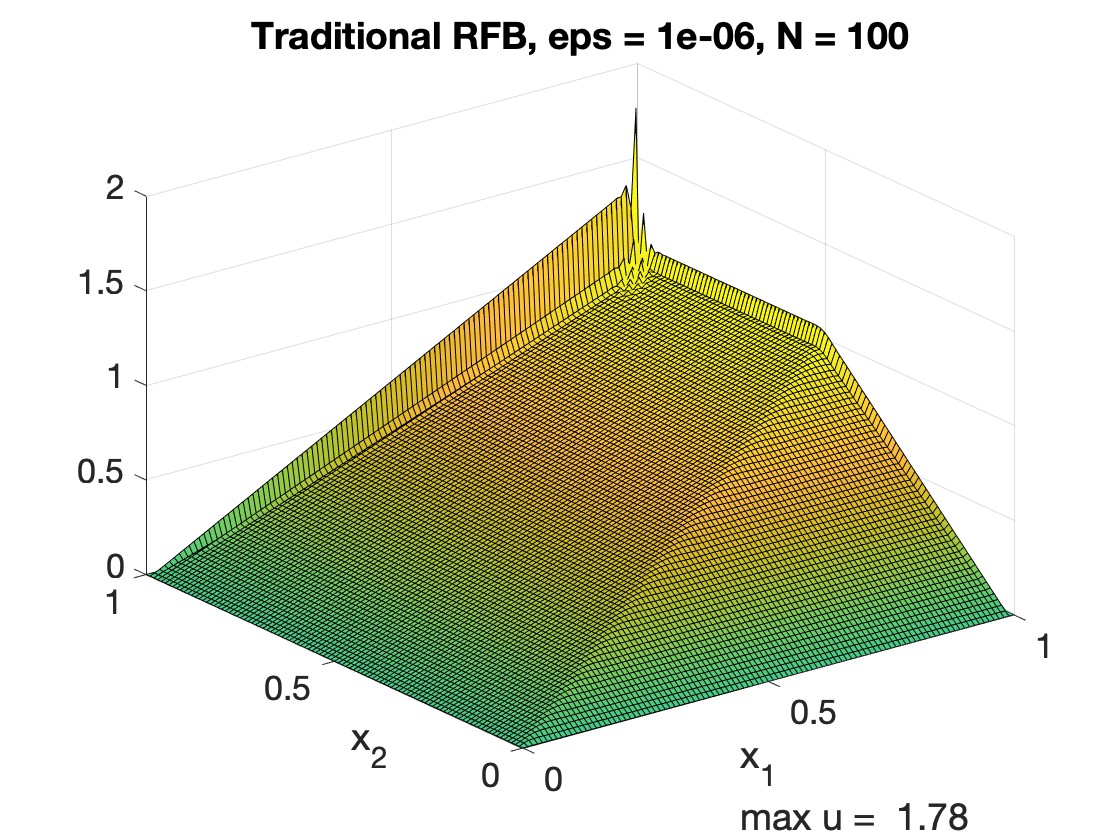}

    \caption{\small\label{F:Brezzi-1} \Change{}{\textit{Example 0,}
    Stabilization with the usual residual-free bubbles.
    Solution of~\eqref{eq:convdiff-stat-strong} on the unit square, with $\epsilon=10^{-6}$, $\aaa=(1, 0.5)$, $c=0$ and $f=1$, using square} meshes for $h=1/50$ (left) and $h=1/100$ (right).
    Even though this solution is much better than the one without any stabilization, we notice still strong oscillations close to the boundary layer, and a maximum value around an 80\% higher than expected.
    }
\end{figure}

\Change{}{
\vspace{0.5 cm}

\textit{Example 0: }
In order to visualize the results obtained with different methods, we consider
Eq.~\eqref{eq:convdiff-stat-strong} on the unit square with $\epsilon = 10^{-6}$,
$\aaa = (1,0.5)$, $c = 0$, $f = 1$, and $g=0$.
Figure~\ref{F:Brezzi-1} shows the solution obtained with the residual-free element
bubble method on square meshes with $h = 1/50$ and $h = 1/100$.
}

Even though this solution is much better than the one without any stabilization, we notice still strong oscillations close to the boundary layer, and the maximum value of $u$ around $1.8$, way above the expected value, which should be below $1$.

We believe that these oscillations appear due to the fact that the residual-free bubbles do not add any \emph{freedom along the edges} of the mesh. In the next section we propose a cure for this, using \emph{patch bubbles}, i.e., 
residual-free bubbles whose supports are unions of two adjacent elements.

\begin{rem}[Space richness\label{rm:SpaceRichness}]\rm
    It is worth noting that in the case of (piecewise) constant data $f$ and $\aaa$, and triangular meshes, the method just presented is equivalent to considering the full $\Hoi[T]$ space instead of the bubble space $B_T=\span \{ \bT \}$. 
    This is due to the fact that for a piecewise linear function $v_L$, we have that $-\epsilon \Delta v_L + \aaa \cdot \nabla  v_L - f$ is constant in each $T$. 

    In order to achieve the same property with quadrilateral meshes, we need four bubbles per element. 
    More precisely, given a quadrilateral element $T$ we let $\varphi_T^i$, $i=1,\dots, 4$ denote the nodal basis functions corresponding to the bilinear elements, and for each $i$, the bubble function $\bT[i]$ is defined as:  
\[
\bT[i] \in \Hoi[T]:\qquad a(\bT[i], \varphi) = (\varphi_T^i, \varphi) \qtfa \varphi \in \Hoi[T].
\]  
Subsequently, the space $B_T = \span\{\bT\}$ is replaced by $\tilde{B}_T = \span\{\bT[i], i=1,\dots, 4\}$. 
The advantage of this modification is that the solution obtained considering the modified space coincides with the solution obtained with $\Hoi[T].$
\end{rem}

\begin{figure}[h!tbp]
    \hfil
    \includegraphics[width=.32\textwidth]{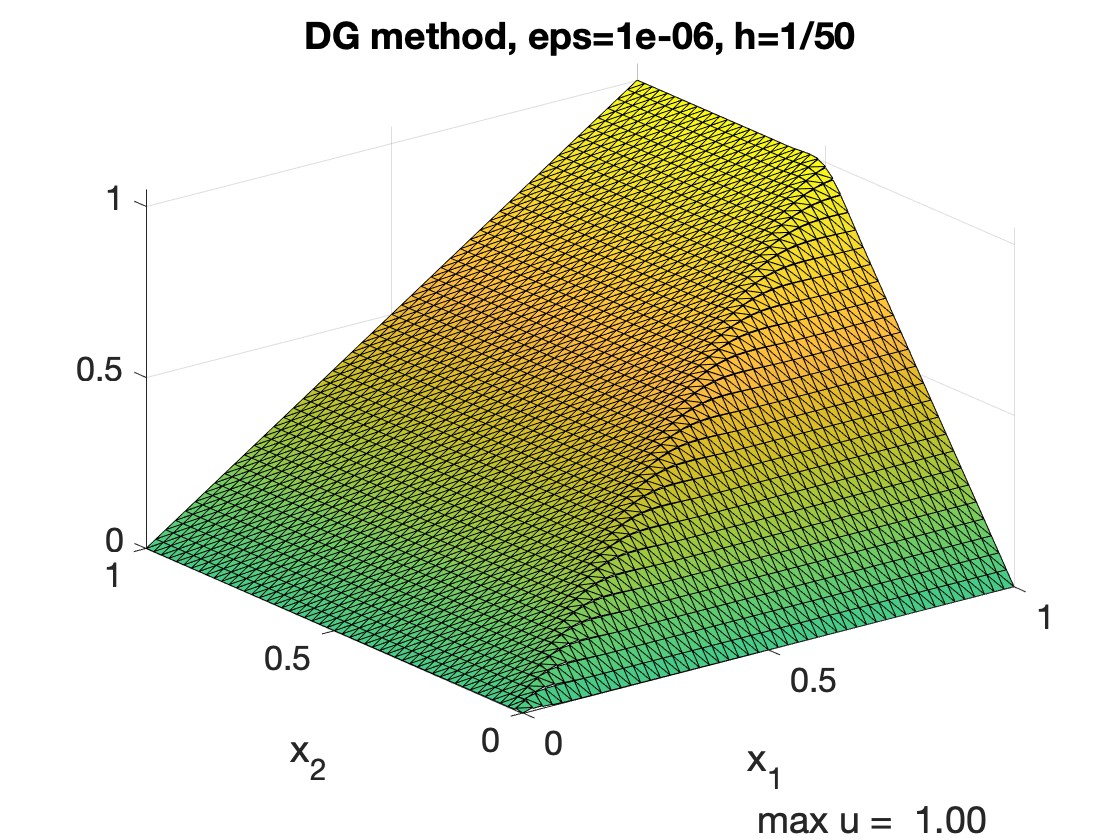}
    \hfil
    \includegraphics[width=.32\textwidth]{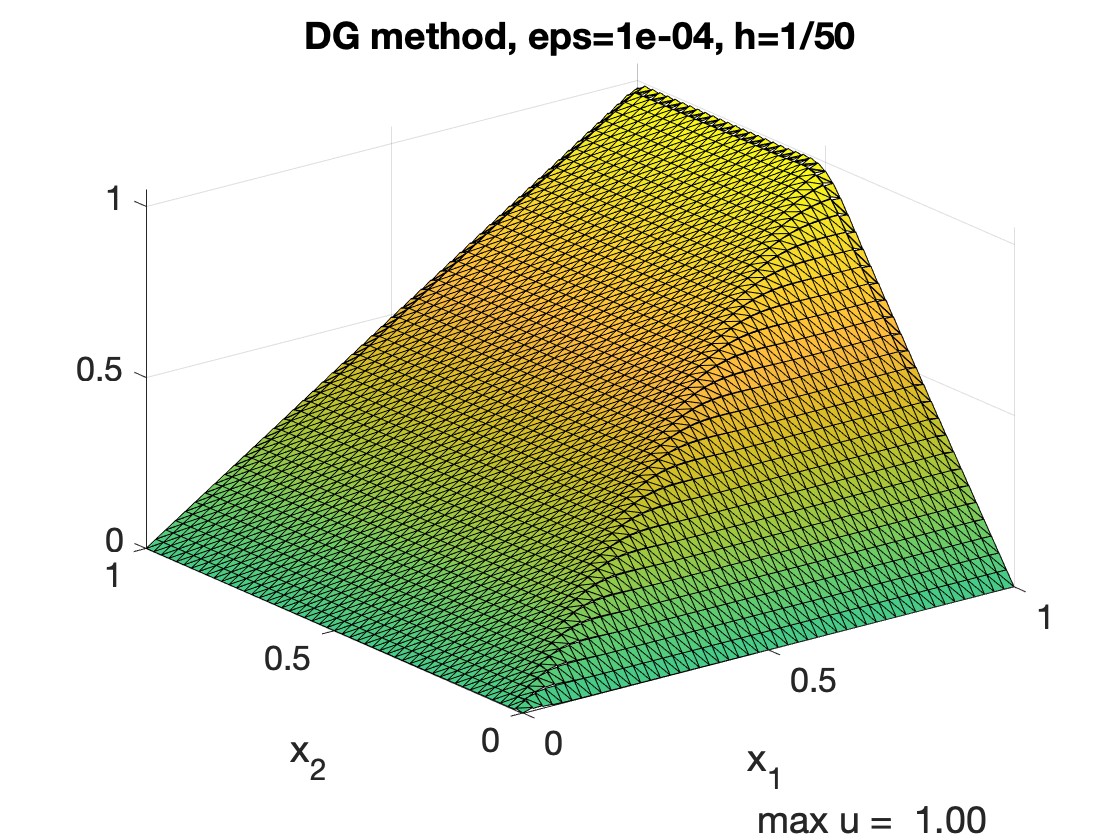}
    \hfil
    \includegraphics[width=.32\textwidth]{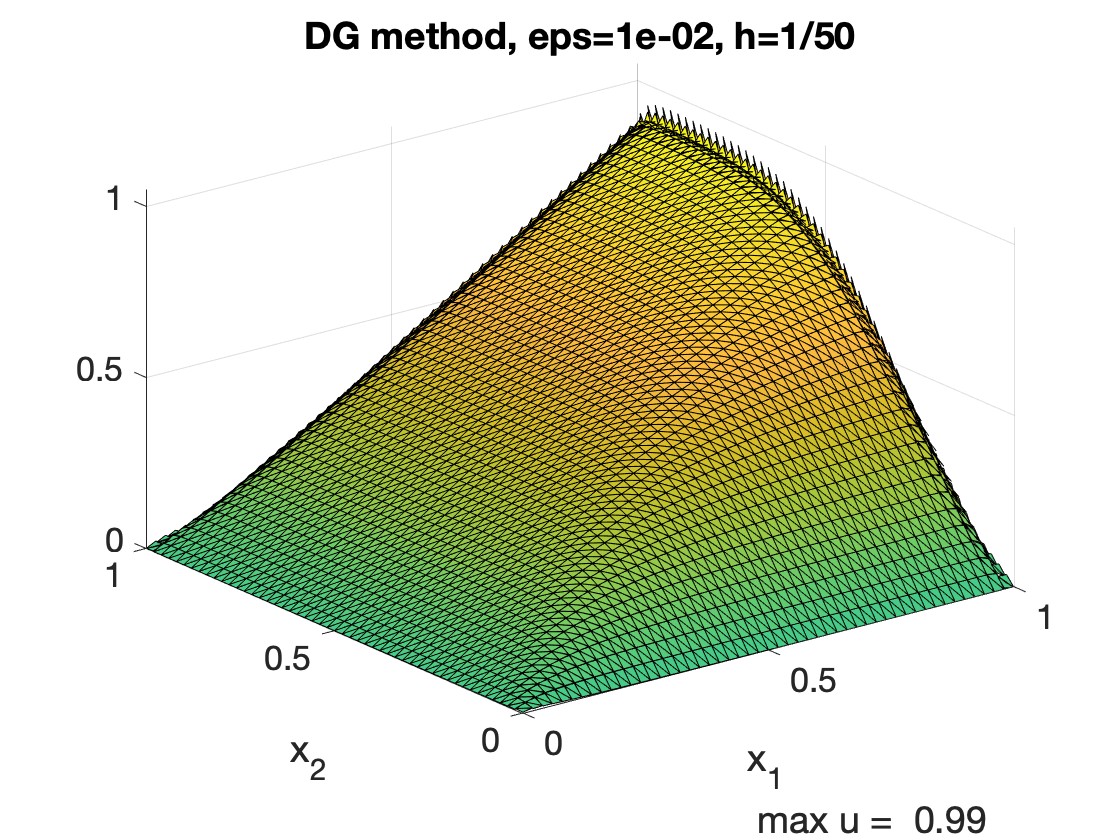}

    \includegraphics[width=.32\textwidth]
{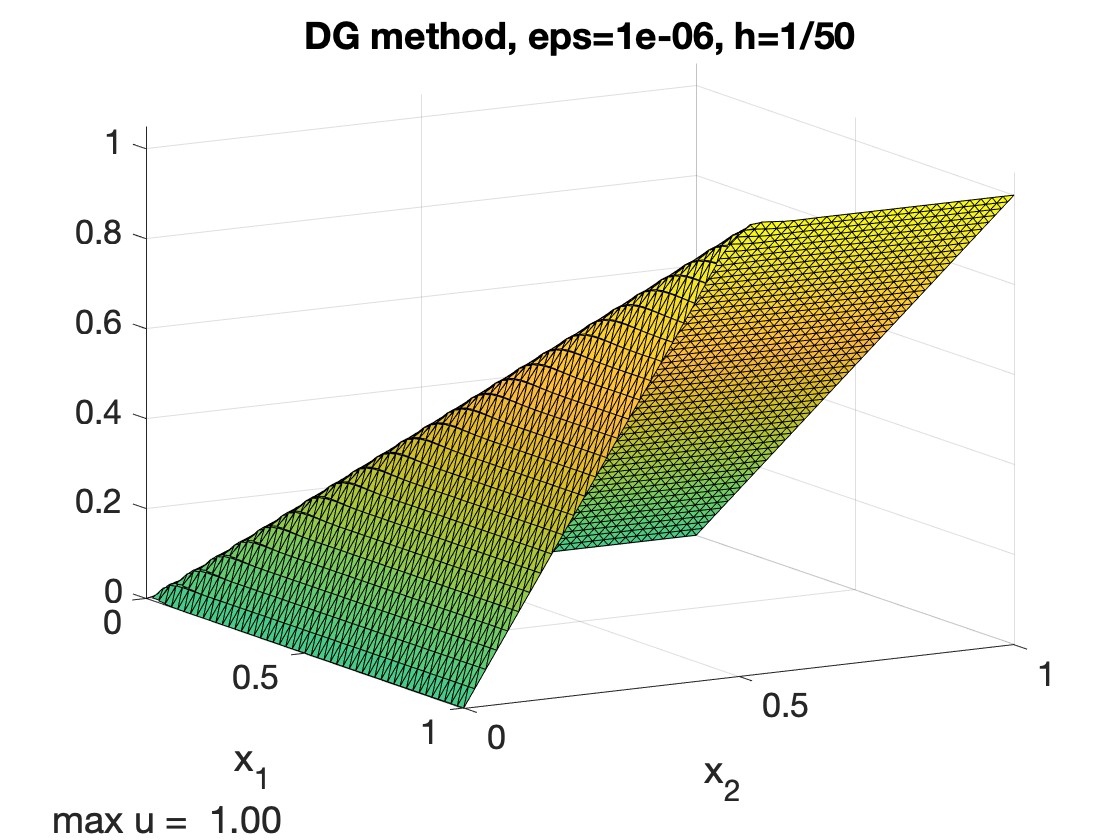}
    \hfil
    \includegraphics[width=.32\textwidth]{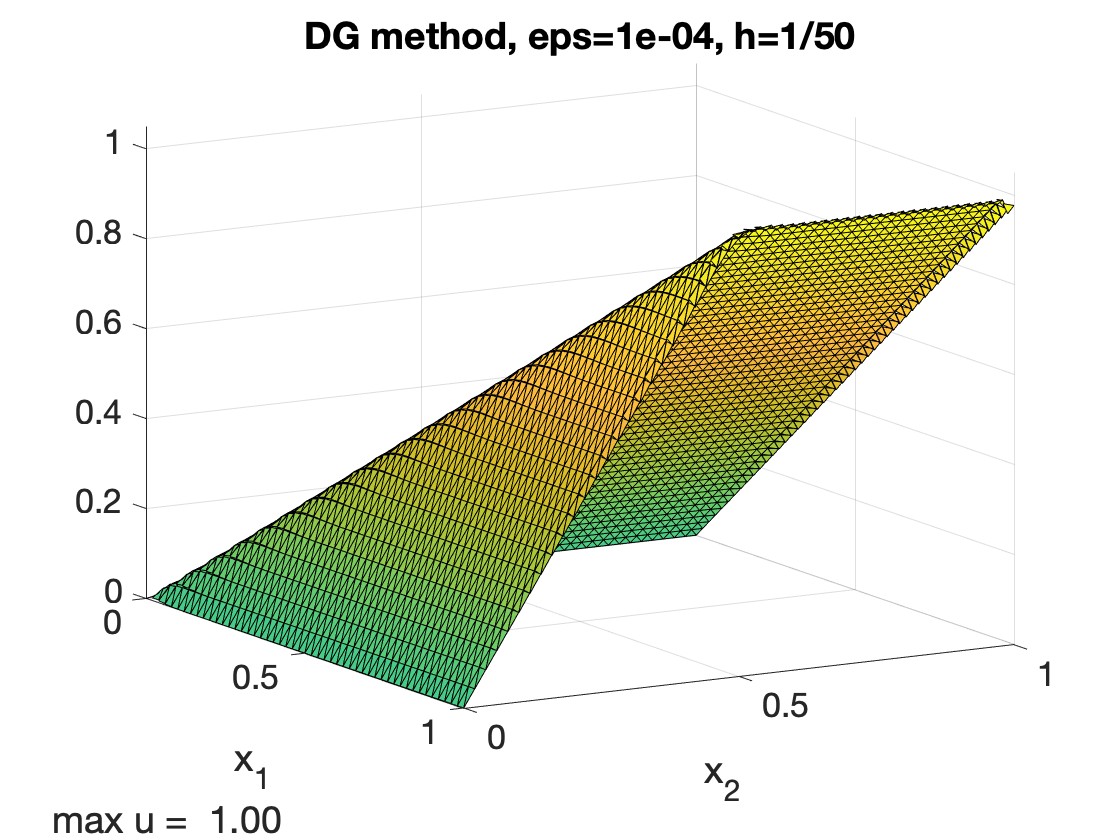}
    \hfil
    \includegraphics[width=.32\textwidth]{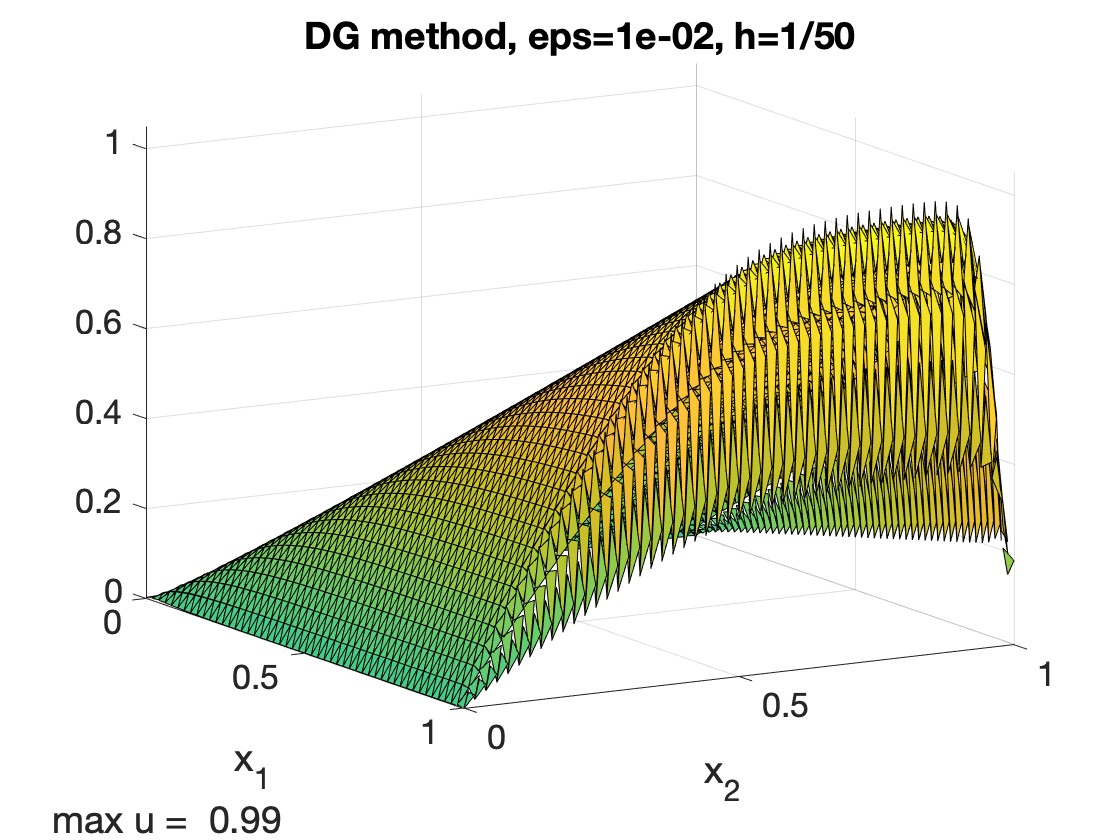}


    \caption{\small\label{F:DG}\Change{}{Example 0,}
    Computation with a DG method. 
    Solution of~{\eqref{eq:convdiff-stat-strong}} in the unit square, with $\aaa=(1, 0.5)$, $c=0$, $f=1$ $g=0$, $\epsilon=10^{-6}$ (left), $\epsilon=10^{-4}$ (middle) and $\epsilon=10^{-2}$ (right), with uniform meshes of size $h=1/50$. Each row corresponds to a different view angle.
    The computation was performed using the toolbox \textsc{festung} {\cite{FESTUNG}}.
    The solutions are excellent for smaller values of $\epsilon$ but some instabilities appear when $\epsilon$ is relatively big.}
\end{figure}

As mentioned in the introduction, another method which exhibits very good performance for advection dominated problems is the Discontinuous Galerkin (DG) method. 
The approximation of the same problem with DG leads to very good solutions for small values of $\epsilon$; see Figure~{\ref{F:DG}}. However, some spurious oscillations appear at the boundary layer for relatively large sizes of $\epsilon$.
The computations were performed using the toolbox \textsc{festung} {\cite{FESTUNG}}, using linear discontinuous finite elements with interior penalization (IIPG), parameter ${\eta=6/h}$, Zalesak Flux limiters with parameters $\gamma=0$, $\beta=10^{-4}$. 
This choice led, in our experience, to the best results.

\section{Residual-Free Patch Bubbles for steady Problems}
\label{S:RFB-patch-stat}

In this section we consider the \delete{following} steady problem~\eqref{eq:convdiff-stat-strong}
\delete{
\begin{equation*} 
    \left\{
\begin{aligned}
    - \epsilon \Delta u + \aaa\cdot \nabla  u  + c\, u &=f \qquad& &\text{in $\Omega$,}\\
        u &= g, \qquad& &\text{on $\partial\Omega$} 
\end{aligned}
    \right.
\end{equation*}
}
with (piecewise) constant $\aaa$, $c$, $f$, and continuous boundary data $g$. 

The associated bilinear form is the same as before:
\begin{equation*}
    a(\psi,\varphi) = \int_\Omega \epsilon \nabla  \psi \cdot \nabla  \varphi + \aaa \cdot \nabla  \psi \, \varphi + c \, \psi \, \varphi,
\end{equation*}
and its weak formulation reads:
\[
\text{Find } u \in \Hgi : \qquad a(u,\varphi) = (f,\varphi) \qtfa \varphi \in \Hoi,
\]
where $(f,\varphi)=\int_\Omega f\varphi$, as usual, 
and $\Hgi = \{ v \in H^1(\Omega) : v_{|\partial \Omega} = g\}$.

\subsection{Bubble Definition}

\Change{In order to simplify the presentation of our main ideas, we focus on $\Omega=(0,1)^2$ and, given $N\in\N$, consider uniform meshes $\T_h$ consisting of $N^2$ square elements of side length $h=1/N$.}{
    We consider a polygonal domain $\Omega \subset \R^2$ which is partitioned into a global mesh $\T_h$ consisting of parallelograms and triangles.
}

We let the finite element space $V_L$ be defined as
\new{
\begin{equation*}
    \begin{aligned}
V_L = \big\{ v \in H^1(\Omega) : \ &v_{|T} \in \Q^1, \text{ for all parallelograms } T \text{ of } \T_h, \\
& v_{|T} \in \PP^1, \text{ for all triangles } T \text{ of } \T_h, \\
&v(x) = g(x), \text{ for all boundary nodes $x$}  \big\},
    \end{aligned}
\end{equation*}
with $\Q^1$ the space of bilinear functions of the form $a_{00} + a_{10} x_1 + a_{01} x_2 + a_{11} x_1 x_2$.
}

Our proposal consists in enriching the discrete space $V_L$ by adding:
\begin{itemize}
\item for each \new{quadrilateral} element $T$, the space $B_T = \span\{\bT[i], i=1,\dots, 4\}$, where $\bT[i]$, $i=1,\dots,4$, are the element bubbles defined in Remark~\ref{rm:SpaceRichness}, i.e., 
\[
\bT[i] \in \Hoi[T]:\qquad a(\bT[i], \varphi) = (\varphi_T^i, \varphi) \qtfa \varphi \in \Hoi[T]
\] 
with $\varphi_T^i$, $i=1,\dots, 4$ the nodal basis functions corresponding to the bilinear elements;
\new{
\item for each triangular element $T$, the space $B_T = \span\{\bT[i], i=1,\dots, 3\}$, where $\bT[i]$, $i=1,\dots,3$, are the element bubbles defined in Remark~\ref{rm:SpaceRichness}, i.e., 
\[
\bT[i] \in \Hoi[T]:\qquad a(\bT[i], \varphi) = (\varphi_T^i, \varphi) \qtfa \varphi \in \Hoi[T]
\] 
with $\varphi_T^i$, $i=1,\dots, 3$ the nodal basis functions corresponding to the linear elements;
}
\item and also for each edge $S$ of the skeleton $\Sigma_h$ of the partition, the space $B_S = \span\{\bS\}$, with the residual-free bubble $\bS$ defined as follows:
\[
\bS\in \Hoi[\patchS]:\qquad a(\bS,\varphi)=(1,\varphi) \qtfa \varphi\in \Hoi[\patchS],
\]
where the patch $\patchS$ is the union of the two elements of the partition sharing the side $S$ (see Figure \ref{F:PatchDomains}). Here, $\Sigma_h$ denotes the set of interior edges of the partition.
\end{itemize}
Summarizing, we set
\[
    \Vrfbn = V_L \oplus V_B \oplus V_S,
    \quad\text{with}\quad
    V_B = \bigoplus_{T\in\T_h} B_T,\tand
    V_S = \bigoplus_{S\in\Sigma_h}B_S.
\]

The discrete problem now reads
\begin{equation}\label{def:u_h}
    \text{Find } u_h \in \Vrfbn:\quad
    a(u_h,v_h) = (f,v_h)\qtfa v_h \in \Vrfbn.    
\end{equation}

The solution $u_h$ can be split as $u_h = u_L + u_B + u_S$, with $u_L \in V_L$, $u_B \in V_B$ and $u_S \in V_S$. Moreover, if 
\begin{equation}\label{def:u_h=uL+B+S}
\begin{aligned}
u_L = \sum_{j=1}^{\#\text{vertices}} u_j \varphi_j,\qquad
u_B = \sum_{j=1}^{\#\text{elements}} \sum_{k=1}^4 u_{T_j}^k \psi_{T_j}^k ,\qquad
u_S = \sum_{j=1}^{\#\text{edges}} u_{S_j} b^{S_j},
\end{aligned}
\end{equation}
the associated linear system can be written as
\begin{equation*}\label{eq: RFBn matricial}
    \begin{bmatrix}
        &\vline&&\vline&\\
        \hspace{0.2 cm}a(\varphi_j,\varphi_i)\hspace{0.2 cm} & \vline & \hspace{0.2 cm}a(\psi_{T_j},\varphi_i)\hspace{0.2 cm} &\vline &\hspace{0.2 cm}a(b^{S_j},\varphi_i)\hspace{0.2 cm}& \\
        &\vline&&\vline&\\
        \hline
        &\vline&&\vline&\\
        \hspace{0.2 cm}a(\varphi_j,\psi_{T_i})\hspace{0.2 cm} & \vline & \hspace{0.2 cm}a(\psi_{T_j},\psi_{T_i})\hspace{0.2 cm}&\vline&\hspace{0.2 cm}a(b^{S_j},\psi_{T_i})\hspace{0.2 cm}&\\
        &\vline&&\vline&\\
        \hline	
        &\vline&&\vline&\\
        \hspace{0.2 cm}a(\varphi_j,b^{S_i})\hspace{0.2 cm} & \vline & \hspace{0.2 cm}a(\psi_{T_j},b^{S_i})\hspace{0.2 cm}&\vline&\hspace{0.2 cm}a(b^{S_j},b^{S_i})\hspace{0.2 cm}&\\
        &\vline&&\vline&\\
    \end{bmatrix}
    \begin{bmatrix}
        \\
        \hspace{0.1 cm}u_j\hspace{0.1 cm} \\
        \\
        \hline
        \\
        \hspace{0.1 cm}u_{T_j}\hspace{0.1 cm}\\
        \\
        \hline
        \\
        \hspace{0.1 cm}u_{S_j}\hspace{0.1 cm}\\
        \\
    \end{bmatrix}
    \hspace{0.0 cm}
    = 
    \hspace{0.0 cm}
    \begin{bmatrix}
        \\	
        \hspace{0.1 cm}(f,\varphi_i)\hspace{0.1 cm} \\
        \\
        \hline
        \\
        \hspace{0.1 cm}(f,\psi_{T_i})\hspace{0.1 cm}\\
        \\
        \hline
        \\
        \hspace{0.1 cm}(f,b^{S_i})\hspace{0.1 cm}\\
        \\
    \end{bmatrix}.
\end{equation*}

As far as we know, the idea of using \emph{patch bubbles} was first introduced by Cangiani and S\"uli 
in~\cite{Cangiani2005long,Cangiani2005short}, in the so called \emph{enhanced residual-free bubbles method} (RFBe),
where they propose to use an exact 1D-solution along the edge $S$, and extend this edge value to a bubble of the patch as a (piecewise) solution to $-\epsilon \Delta u + \aaa \cdot \nabla  u=0$ 
with homogeneous boundary values.
They also prove some error estimates and highlight the power of using \emph{enhanced} bubbles, besides showing some impressive numerical experiments.

Patch bubbles have a disadvantage compared to element bubbles though: their incorporation requires applying 
the full bilinear form to various combinations of functions. Specifically, it is necessary to compute the interactions 
of patch bubbles with themselves, with element bubbles, and with standard finite element basis functions. This involves 
evaluating the bilinear operator associated with the weak form, i.e., computing $a(f_i,f_j)$ where the functions 
$f_i$ and $f_j$ may represent affine functions $\varphi$, element bubbles $\bT$, or patch bubbles $\bS$; we call these integrals \emph{stabilization terms}.

A particular challenge arises in some of these interactions, because integration by parts cannot always simplify the computations to involve only the integral of the bubbles, as is the case in~\eqref{eq:bT-identity-1}--\eqref{eq:bT-identity-2}. To illustrate, consider the integral in the bilinear form \eqref{eq:bilinear-op-a} associated with the second-order term: ${\epsilon\int\nabla  u \cdot \nabla  v}$.
In the interaction between a patch bubble $\bS$ and a nodal basis function $\varphi$, the following expressions are obtained:

\begin{equation*}
\begin{aligned}
 \int_{\patchS} \epsilon\nabla  \bS \cdot \nabla  \varphi &= \int_{\partial\patchS}\epsilon (\nabla  \bS \cdot \nnn) \varphi -\int_{\patchS}\epsilon(\varphi \Delta \bS), \\
    \int_{\patchS}\epsilon\nabla  \varphi \cdot \nabla  \bS &= \sum_{T \subset\patchS} \int_{S}\epsilon \bS(\nabla \varphi \cdot \nnn).
\end{aligned}
\end{equation*}

Therefore, we face the problem of computing the bubbles with more precision. However, the bubbles are solutions to a problem as difficult as the original one; albeit with an important difference: the domain is now of size $h$.
Cangiani and S\"uli use Shishkin meshes in order to compute these bubbles.

Our approach, instead, leverages \emph{recursion}. The values of $a(f_i,f_j)$ are computed within the function called $\texttt{element\_contribs}$, using the procedure described in Algorithm \ref{Alg:compute bubbles} 
, referred to as $\texttt{stab\_terms}$. More details are provided in the following section.

\subsection{Recursive Bubble Computation}

As we mentioned in the previous section, to compute the patch bubbles we propose a recursive algorithm. More precisely, the algorithm proceeds as follows.

\Change{When we need to compute the bubbles, we generate finite element partitions on the \emph{patches}, and use the same 
method with patch bubbles to compute the bubbles on the global mesh. The number of different patches is finite, and much smaller than the number of elements on regular meshes, and can be done in parallel, or \emph{offline}.}{
    In order to compute the bubbles, uniform meshes are generated on each element of the patch, by refining all the sides into $M$ elements of equal length, and splitting the elements with lines parallel to the original sides. This procedure, which we call \texttt{refine\_patch(patch,$M$)} below, creates elements that are all similar to the original elements of the mesh. Continuing in the recursion, all patches will be either similar to the original patches, or will contain two elements that are both similar to one element of the original patch; with three possible orientations for triangles and two posible orientations for parallelograms. When using square meshes, there will be only two classes of patches (see Section~\ref{S:implementation} below). When using simplicial meshes, there will be seven kinds of patches for each pair of elements (see Figure~\ref{F:PatchDomains})
}
This leads to a recursive algorithm, which will reach the limit of \Change{$h < \epsilon/|\aaa|$}{$Pe_h<1$} in $O(\hspace{0.05 cm} \log (\hspace{0.05 cm} \Change{|\aaa| h \hspace{0.05 cm}/\hspace{0.05 cm}\epsilon}{\Pe_h} ))$ recursive steps.

\begin{figure}
    \centering
    \begin{tikzpicture}[scale=0.7]

        \node[font=\footnotesize] at (0,-0.5) {$$};
      \node at (1,1.2) {\textcolor{blue}{$\Omega$}};
        
\coordinate (A) at (0,0);
\coordinate (B) at (4,0);
\coordinate (C) at (6,1);
\coordinate (D) at (4,2);
\draw[blue, very thick]
  (A) -- (B) -- (C) -- (D) -- cycle;

\draw[blue, very thick]
  (4,0) -- (4,2);
  \draw[fill=blue,opacity=0.4,thick, blue] (1,0) rectangle (2,0.5);
\draw[fill=black,opacity=0.4, thick]
   (2,0) -- (3,0.5) -- (3,1) -- (2,0.5) -- cycle;
\draw[fill=orange,opacity=0.4, thick]
   (2,1) -- (3,1) -- (4,1.5) -- (3,1.5) -- cycle;
\draw[fill=green,opacity=0.4, thick]
   (3,1) -- (4,1) -- (4.5,1.25) -- (4,1.5) -- cycle;
\draw[fill=blue,opacity=0.4, thick]
   (4,1) -- (4.5,1.25) -- (5,1) -- (4.5,0.75) -- cycle;
\draw[fill=orange,opacity=0.4, thick]
   (4.5,0.75) -- (5,1) -- (5,.5) -- (4.5,0.25) -- cycle;
\draw[fill=black,opacity=0.4, thick]
   (5,1) -- (5,1.5) -- (5.5,1.25) -- (5.5,0.75) -- cycle;

    \draw[black!70!, thick](1,0) -- (1,0.5);
    \draw[black!70!, thick](2,0) -- (2,1);
    \draw[black!70!, thick](3,0) -- (3,1.5);

    \draw[black!70!, thick](1,0) -- (4,1.5);
    \draw[black!70!, thick](2,0) -- (4,1);
    \draw[black!70!, thick](3,0) -- (4,0.5);

    \draw[black!70!, thick](1,0.5) -- (4,0.5);
    \draw[black!70!, thick](2,1.0) -- (4,1.0) ;
    \draw[black!70!, thick](3,1.5) -- (4,1.5);

    \draw[black!70!, thick](4,0.5) -- (4.5,0.25);
    \draw[black!70!, thick](4,1.0) -- (5.0,0.5);
    \draw[black!70!, thick](4,1.5) -- (5.5,0.75);

    \draw[black!70!, thick](4.5,0.25) -- (4.5,1.75);
    \draw[black!70!, thick](5.0,0.5)  -- (5.0,1.5);
    \draw[black!70!, thick](5.5,0.75) -- (5.5,1.25);

    \draw[black!70!, thick](4,0.5) -- (5.5,1.25);
    \draw[black!70!, thick](4,1.0) -- (5.0,1.5);
    \draw[black!70!, thick](4,1.5) -- (4.5,1.75);
\end{tikzpicture}
\begin{tikzpicture}[scale=0.3]
\node at (6,3) {\textcolor{blue}{$\Omega$}};

\coordinate (A) at (0,0);
\coordinate (B) at (-1,4);
\coordinate (C) at (-5,3);
\coordinate (D) at (-4,-1);

\coordinate (E) at (8,-2);
\coordinate (F) at (7,2);
\draw[blue, very thick]
  (A) -- (B) -- (C) -- (D) -- cycle;

\draw[blue, very thick]
  (A) -- (B) -- (F) -- (E) -- cycle;

\draw[fill=orange,opacity=0.4, thick]
    (-4,3.25) -- (-3.5,1.25) -- (-2.5,1.5) -- (-3,3.5) --cycle;
\draw[fill=blue,opacity=0.4, thick]
    (-3.5,1.25) -- (-3.25,0.25) -- (-1.25,0.75) -- (-1.5,1.75) --cycle;
\draw[fill=green,opacity=0.4, thick]
    (-1.5,1.75) -- (-0.5,2) -- (1.5,1.5) -- (1.25,2.5) --
    (-0.75,3) -- (-1.75,2.75) -- cycle;
\draw[fill=blue,opacity=0.4, thick]
    (-0.5,2) -- (-0.25,1) -- (3.75,0) -- (3.5,1) --cycle;
    \draw[fill=orange,opacity=0.4, thick]
    (3.75,0) -- (5.75,-0.5) -- (5.25,1.5) -- (3.25,2) --cycle;

    \foreach \i in {1,...,3} {
  \draw[black!70!, thick] (-\i,-0.25*\i) -- ++(-1,4);
  }

    \foreach \i in {1,...,3} {
  \draw[black!70!, thick] (-0.25*\i,\i) -- ++(-4,-1);
  }

    \foreach \i in {1,...,3} {
  \draw[black!70!, thick] (2*\i,-0.5*\i) -- ++(-1,4);
  }

    \foreach \i in {1,...,3} {
  \draw[black!70!, thick] (-0.25*\i,\i) -- ++(8,-2);
  }

\end{tikzpicture}
\begin{tikzpicture}[scale=0.4]
\node at (-2,0) {\textcolor{blue}{$\Omega$}};

\coordinate (A) at (0,0);
\coordinate (B) at (-1,4);
\coordinate (C) at (-4,2);

\coordinate (D) at (4,-1);
\coordinate (E) at (3,3);
\draw[blue, very thick]
  (A) -- (B) -- (C) -- cycle;

\draw[blue, very thick]
  (A) -- (B) -- (E) -- (D) -- cycle;

  \draw[fill=orange,opacity=0.4, thick]
    (-13/4,2.5) -- (-10/4,3) -- (-9/4,2) -- (-3,1.5) -- cycle;
\draw[fill=black,opacity=0.4, thick]
    (-9/4,2) -- (-2,1) -- (-1,0.5) -- (-1.25,1.5) -- cycle;
\draw[fill=blue,opacity=0.4, thick]
    (-10/4,3) -- (-6/4,2.5) -- (-0.75,3) -- (-7/4,3.5) -- cycle;
\draw[fill=green,opacity=0.4, thick]
    (-6/4,2.5) -- (-0.75,3) -- (0.25,2.75) -- (0.5,1.75) -- (-0.5,2) -- cycle;
\draw[fill=blue,opacity=0.4, thick]
    (-0.5,2) -- (-0.25,1) -- (1.75,0.5) -- (1.5,1.5) -- cycle;
\draw[fill=orange,opacity=0.4, thick]
    (1.75,0.5) -- (1.25,2.5) -- (2.25,2.25) -- (2.75,.25) --cycle;
    \draw[black!70!, thick](-1,0.5) -- (-1/4,1);
    \draw[black!70!, thick](-2,1) -- (-2/4,2);
    \draw[black!70!, thick](-3,1.5) -- (-3/4,3);

    \draw[black!70!, thick](-1,0.5) -- (-7/4,3.5);
    \draw[black!70!, thick](-2,1) -- (-10/4,3);
    \draw[black!70!, thick](-3,1.5) -- (-13/4,2.5);

    \draw[black!70!, thick](-7/4,3.5)  -- (-3/4,3);
    \draw[black!70!, thick](-10/4,3)   -- (-2/4,2);
    \draw[black!70!, thick](-13/4,2.5) -- (-1/4,1);

    \foreach \i in {1,...,3} {
  \draw[black!70!, thick] (-0.25*\i,\i) -- ++(4,-1);
  }

    \foreach \i in {1,...,3} {
  \draw[black!70!, thick] (\i,-0.25*\i) -- ++(-1,4);
  }
\end{tikzpicture}
\caption{\new{Examples of patch domains. Seven possible types of patches will appear for each patch of the \Change{coarsest}{global} mesh when using simplices (left). One similar to the original one, with the interelement edge contained in the interelement edge of the original patch.
    Three composed of two elements similar to one of the original elements, and three composed of two elements similar to the other original element. These only differ in the orientation. In the case of quadrilaterals, only five types of patches will appear (middle); and six in the case of mixed patches (right).
    }}
    \label{F:PatchDomains}
\end{figure}

In the experiments below,\new{
    in order to make a first evaluation of the performance of this method,
    }
we consider uniform meshes of squares in $\Omega = (0,1)^2$, thereby obtaining only two reference patches.
\delete{If we were considering a different partition with simplices and quadrilaterals, there would be more reference patches, but inside the recursion, for each element, all the reference sub-patches would be equal.}

The algorithm, which from now on we refer to as the BMZ (Bubble-Mesh-Zoom) method, is based on the idea of zooming in locally to solve bubble subproblems over each mesh element or patch. It is defined through two functions, which we describe in Algorithms~\ref{Alg:stabilized fem} and~\ref{Alg:compute bubbles}. 

\new{
\begin{algorithm}[H]
    \small
    \caption{\small\texttt{BMZ}.}
    \label{Alg:stabilized fem}
    \SetAlgoLined
    \KwIn{$\epsilon$, $\aaa$, $c$, $\Omega$, \new{\texttt{global\_mesh}, }$h$, $f$, $g$, M.}
    \KwOut{Computed solution: $u$. }

    \uIf{$\Change{h < {\epsilon}/{|\mathbf{a}|}}{\Pe_h <1}$}{
        Compute $u$ using the standard Galerkin method\;
    }
    \Else{
        \new{
        $stab\_terms \gets \texttt{element\_contribs}(\epsilon, \aaa, c, \texttt{global\_mesh},M)$\;
        Add  $stab\_terms$ to global matrix\;

        }
        Compute $u$ using the stabilized method with $stab\_terms$\;
    }
\end{algorithm}
}

\bigskip

\new{
\begin{algorithm}[H]
    \small
    \caption{\small\texttt{element\_contribs}.}
    \label{Alg:compute bubbles}
    \SetAlgoLined
    \textbf{Function} \texttt{element\_contribs}($\epsilon, \mathbf{a}, c, \texttt{mesh}, M$)

    \KwOut{Stabilization terms: $stab\_terms$. }
    
    
    \For{patch in \{classes of patches\}}{
        $h$ $\gets$ size of \texttt{patch}\; 
        
        \texttt{patch\_submesh} $\gets$ \texttt{refine\_patch}(\texttt{patch}, $M$);
        
            $h_{\text{new}} \gets {h}/{M}$\;
        
    \uIf{$\Change{h_{\text{new}} < {\epsilon}/{|\mathbf{a}|}}{Pe_{h_\text{new}}<1}$}{
        Compute the bubbles $\bT$, $\bS$ using the standard Galerkin method\;
        }
                \Else{
                        $stab\_terms$$\gets$\texttt{element\_contribs}$(\epsilon, \mathbf{a}, c, \texttt{patch\_submesh},M)$\;
        
        Compute the bubbles using the stabilized method with $stab\_terms$\;
    }}
    Calculate \texttt{stab\_terms}:
    $a(\bT, \bS)$, $a(\bS, \bT)$, $a(\bS, \bS)$, $\int\bT$, $\int\bS$, \dots\;
    \Return \texttt{stab\_terms}\;

\end{algorithm}
}

\bigskip

\begin{figure}[h!tbp]
    
    \hfil
    \includegraphics[width=.32\textwidth]{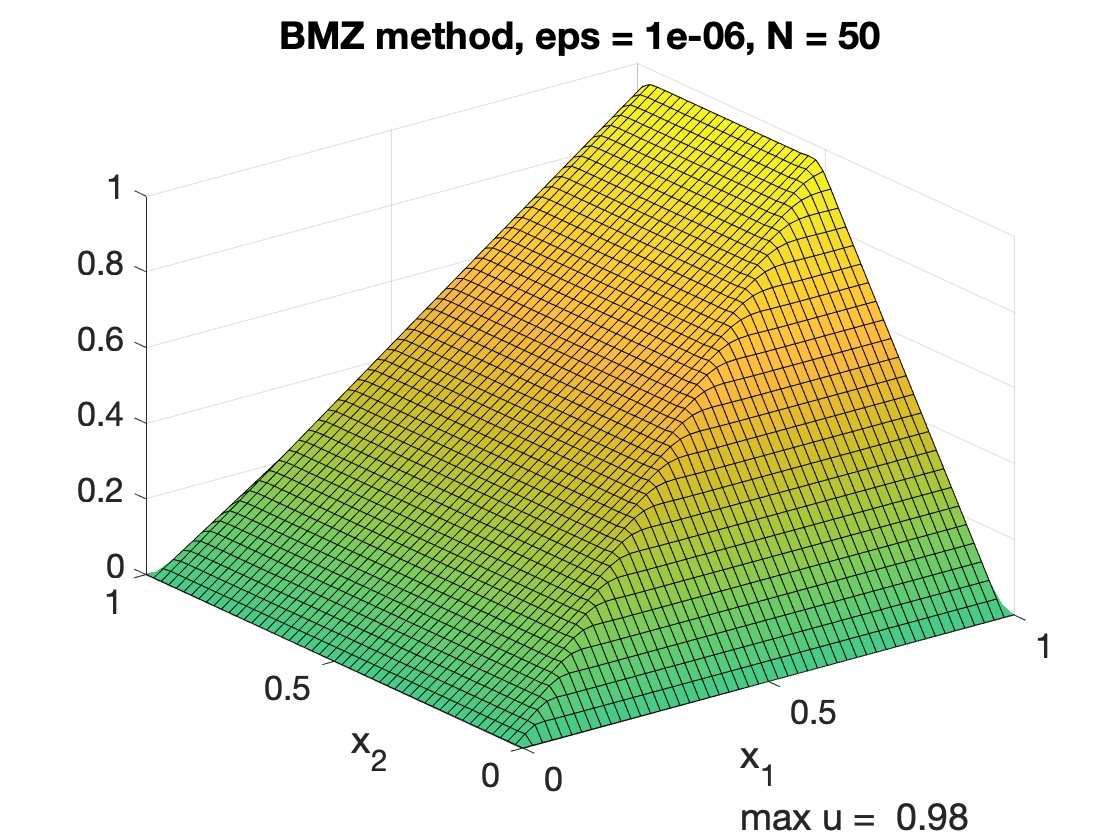}
    \hfil
    \includegraphics[width=.32\textwidth]{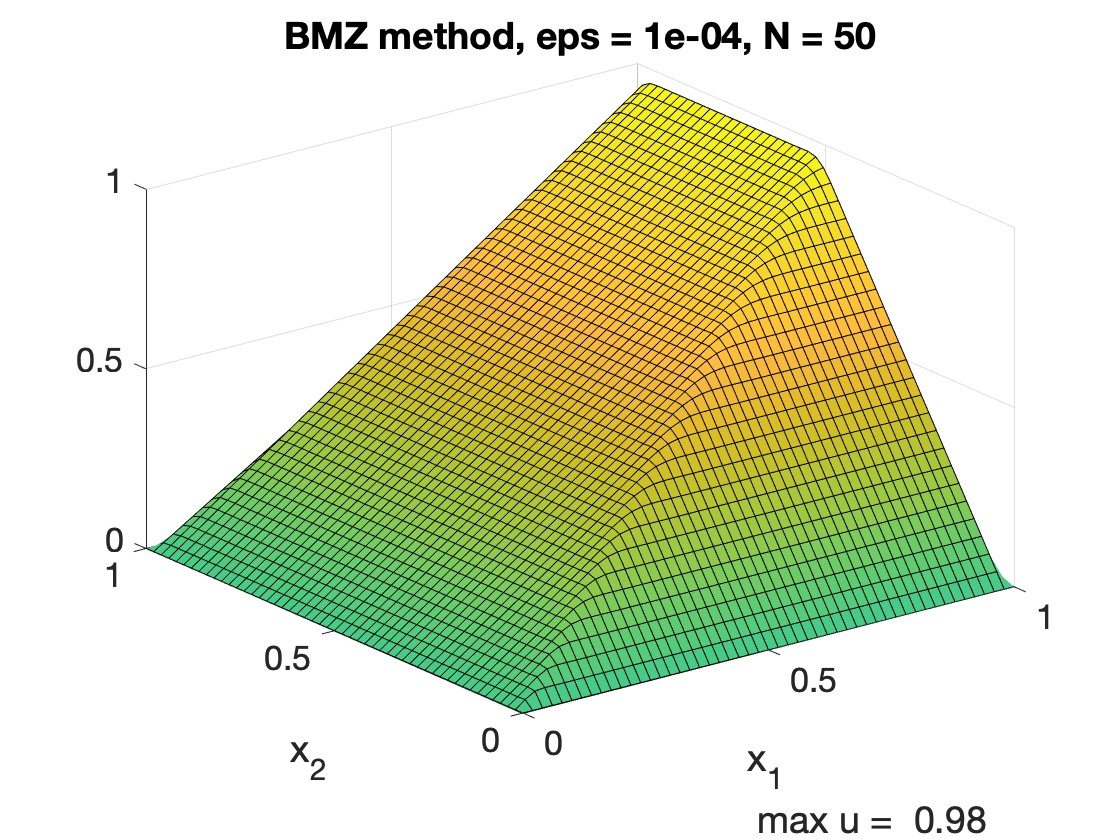}
        \hfil
    \includegraphics[width=.32\textwidth]{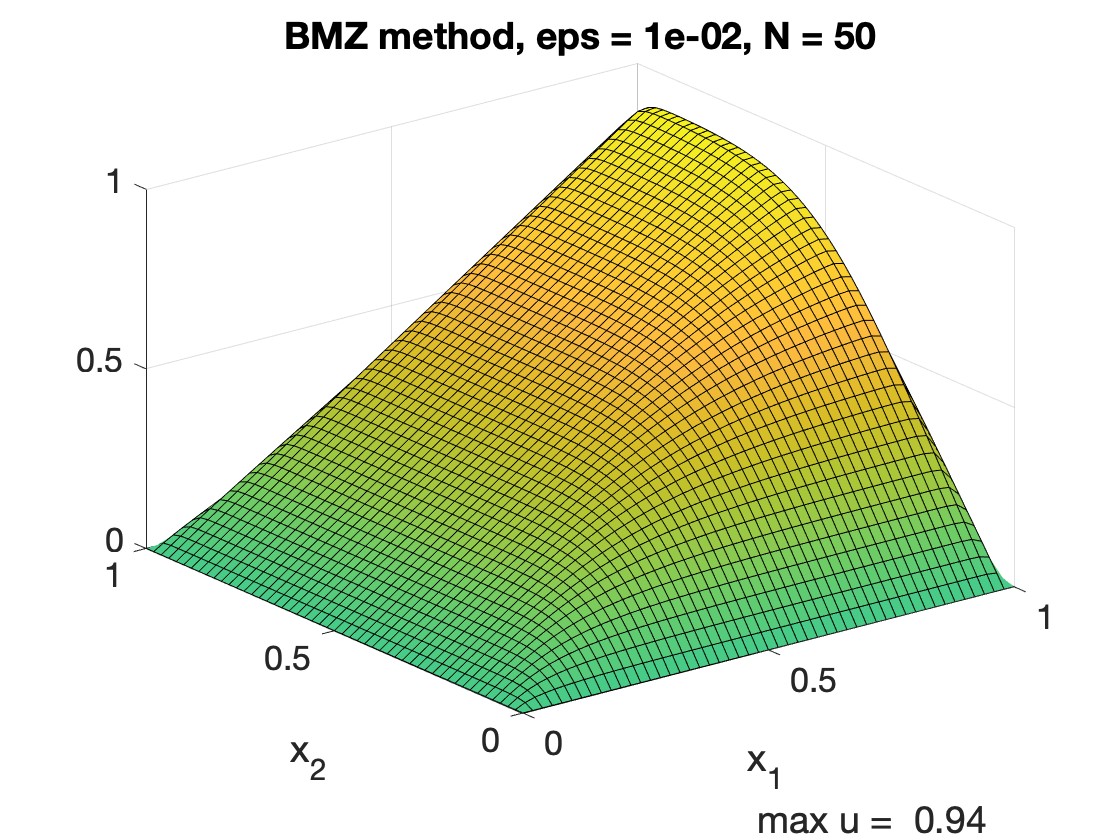}

    \includegraphics[width=.32\textwidth]{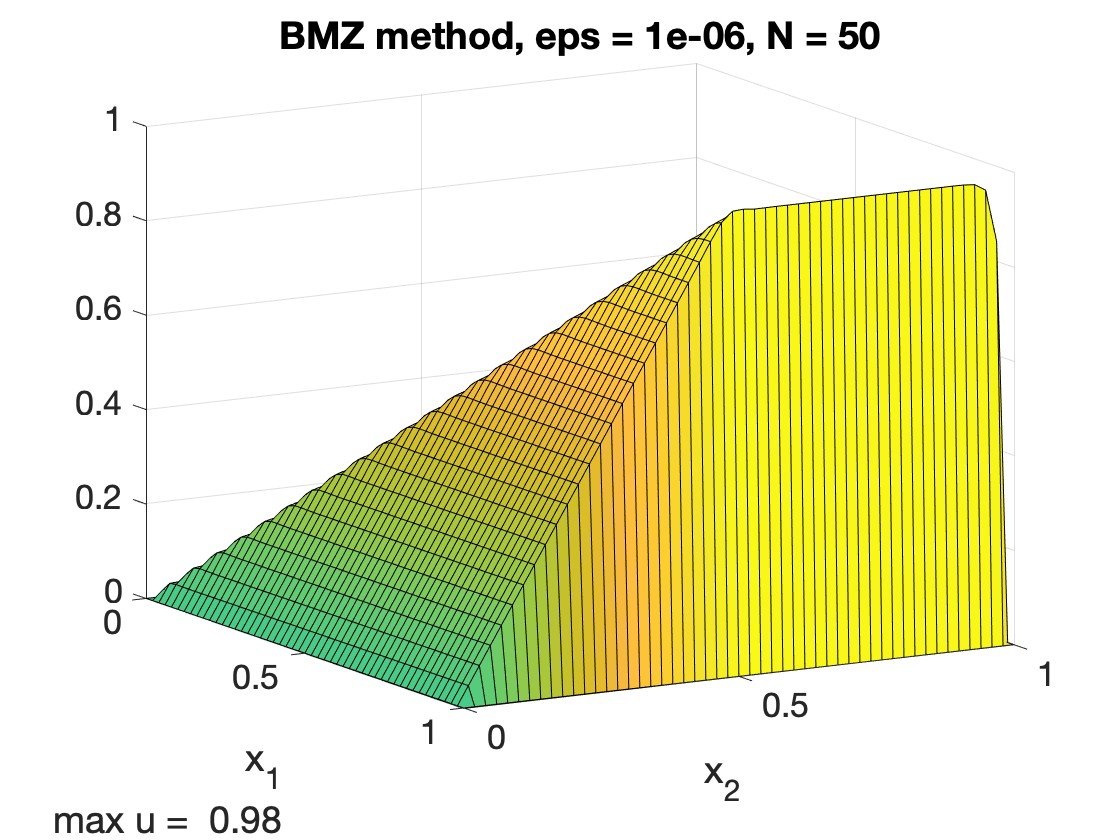}
    \hfil
    \includegraphics[width=.32\textwidth]{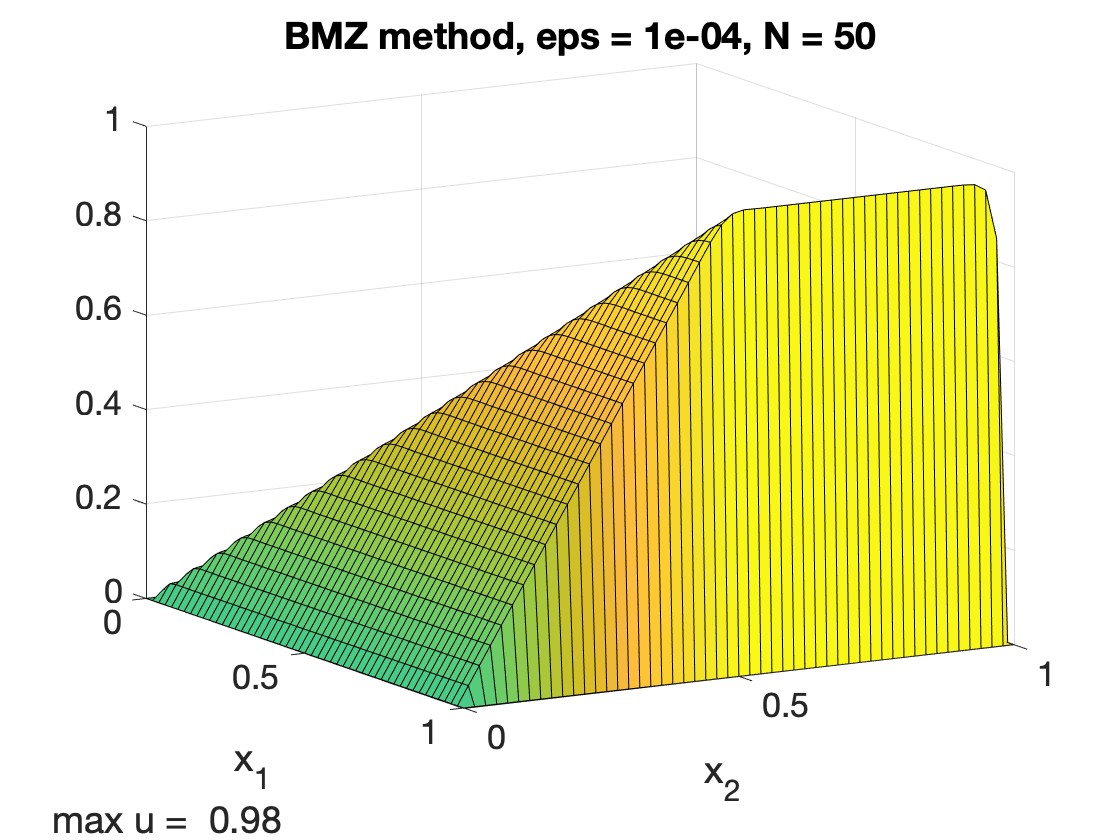}
        \hfil
    \includegraphics[width=.32\textwidth]{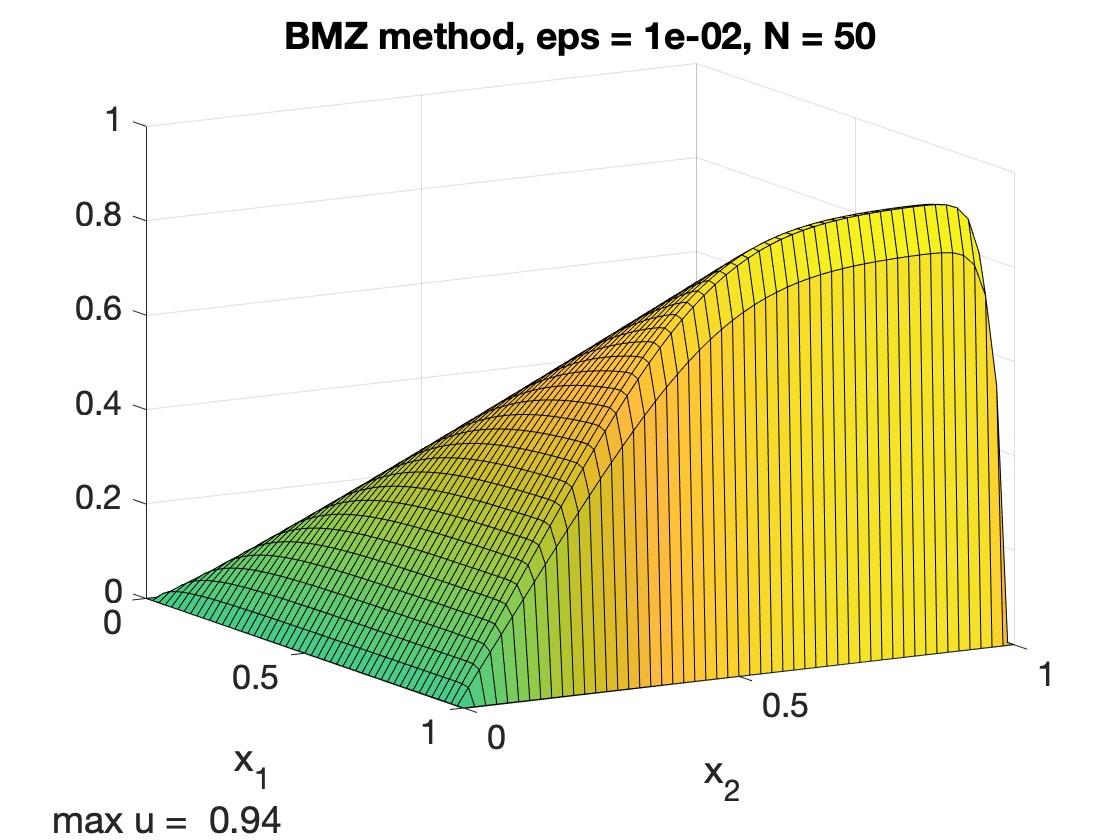}

    \caption{\small\label{F:BMZ-1}\Change{}{\textit{Example 0,}}
    Stabilization with the patch residual-free bubbles.
    Solution of~{\eqref{eq:convdiff-stat-strong}} on the unit square, with $\aaa=(1, 0.5)$, $\epsilon=10^{-6}$ (left), $\epsilon=10^{-4}$ (middle) and $\epsilon=10^{-2}$ (right), with uniform square meshes of size $h=1/50$. Each row corresponds to a different view angle.
    This solution looks strikingly good, with no oscillations whatsoever.
    }
\end{figure}

The \Change{behavior of}{solutions obtained with} these algorithms, \Change{}{for Example 0,} can be observed in Figure~\ref{F:BMZ-1}, \Change{where we show the solution of~\eqref{eq:convdiff-stat-strong} on the unit square, with $\epsilon=10^{-6}$ and $\aaa=(1, 0.5)$ on quadrilateral}{using square} meshes for $h=1/50$ and $h=1/100$; see Section~\ref{S:implementation} for implementation details \Change{}{and the bubbles obtained in the intermediate steps of the BMZ algorithm}.
This solution looks strikingly good, with no oscillations whatsoever, and the maximum value of $u$ approximately $0.9821$ ($h=1/50$) and $0.9916$ ($h=1/100$). It is also noticeable the minimal to no smearing in the solutions.
Very similar results are obtained with the method by Cangiani and S\"uli~\cite{Cangiani2005long,Cangiani2005short}.

\begin{rem}[Non-constant flow fields]
    \label{rem:Non-constant flow fields}
\rm  The algorithm is also capable of coping with 
spatially varying \new{smooth} velocity fields\delete{ in general}. In such cases, Algorithm~\eqref{Alg:stabilized fem} computes the mean value 
of \(\aaa\) within each element and uses it as input for the function \texttt{element\_contribs} (Algorithm~\eqref{Alg:compute bubbles}). \Change{}{Specifically, lines 4--5 of Algorithm~\ref{Alg:stabilized fem} are replaced by the following lines of code:}

\new{
\begin{algorithm*}[H]
    \small
    \SetAlgoLined
    \For{element in \texttt{global\_mesh}}{
    
        $\bar{\aaa}$ $\gets$ mean of \aaa\;
    \For{patch containing element}{
        $stab\_terms$ $\gets$ \texttt{element\_contribs}$(\epsilon, \bar{\aaa}, c,patch,M)$\;
        
        Add  $stab\_terms$ to global matrix\;
    }}
\end{algorithm*}}
For improved performance, the computation of the stabilization terms \Change{can}{could} be parallelized across different cpu's, enabling a more efficient implementation. \Change{}{This would require a careful implementation to avoid imbalance in the computational load; this aspect falls beyond the scope of this article.}

\new{Taking advantage of the zooming procedure to better resolve the velocities with fine scales would be possible, and an interesting idea to explore in future research. In particular, it would be very interesting to consider advection-dominated equations where the velocity field already has a steep boundary layer.}
\end{rem}

\bigskip

\subsection{Numerical Experiments}

    The goal of this section is to illustrate the performance of this idea by solving some steady problems and comparing with the other stabilized methods using bubbles.
We solve on a square with uniform meshes in order to assess the potential power of this method and leave the implementation to more general finite elements, on arbitrary polygonal domains for a future article.

\subsubsection{Implementation on squares}
\label{S:implementation}


For the numerical implementation of the proposed method we consider $\Omega$ as the  square $(0,L)\times(0,L)$, usually with $L=1$ . The partition consists of squares with sides parallel to the axes. For the finite element space, we use bilinear functions. Based on the definition of the space \Vrfbn, we include bubbles of type $\bT$ with support on a single element 
and bubbles $\bS$ with support on the union of two adjacent elements. Consequently, the development of the code requires solving problem~\eqref{eq:convdiff-stat-strong} with 
$f\equiv 1$ and $g\equiv 0$ in the  domains $Q,\patchS[1],\patchS[2]$,
where $Q := (0,h)\times(0,h)$, $\patchS[1] := (0,2h)\times(0,h)$ and $\patchS[2] := (0,h)\times (0,2h).$
\begin{wrapfigure}{r}{0.4\textwidth}
    \centering
    \vspace{- 0.0 cm}	
    \begin{tikzpicture}[scale=2.7]
        \foreach \x in {0,0.2,...,1}
        \draw (\x,0) -- (\x,1);
        \foreach \y in {0,0.2,...,1}
        \draw (0,\y) -- (1,\y);
        
        \draw[thick, blue] (0,0) rectangle (1,1) node[above right] {$\Omega$};
        
        \fill[green!50!black, opacity=0.2] (0.01,0.01) rectangle (0.39,0.19);
        \draw[thick,green!50!black] (0.01,0.01) rectangle (0.39,0.19)node[ below right=-1.5pt] {$\omega_S^1$};
        
        \fill[orange, opacity=0.2] (0.01,0.01) rectangle (0.19,0.39);
        \draw[thick,orange] (0.01,0.01) rectangle (0.19,0.39)node[above left = -2pt] {$\omega_S^2$};
        
        \fill[purple, opacity=0.2] (0.015,0.015) rectangle (0.185,0.185);		
        \draw[line width=0.5pt, purple] (0.015,0.015) rectangle (0.185,0.185) node[midway] {$Q$};		
        
        \node[font=\footnotesize] at (1,-0.06) {$1$};
        
        \node[font=\footnotesize] at (-0.06,1) {$1$};
        
        \node[font=\footnotesize] at (-0.03,-0.06) {$0$};	
    \end{tikzpicture}

\caption{\small\label{F:Ref-Domains} Reference domains.}
\vspace{- 0.5 cm}
\end{wrapfigure}

The bubbles are generated 
by solving a system of order $M^2\times M^2$ for each of the four element bubbles and also for each of the two edge bubbles. 
Here, $M$ is a fixed small value; in the experiments presented in this work, $M=20$ is used. 

The main idea can be summarized as follows: to compute the bubbles, the domain---originally of size $h$---is subdivided (uniformly) into elements of size $h/M$ .

If $h/M$ is less than $\epsilon/|\aaa |$,\Change{}{ that is ${Pe_{h/M} <1}$}, the bubbles  $\bT$ and $\bS$ are obtained by solving the problem using the standard Galerkin method on the generated mesh, and the interactions between these functions and the bilinear basis functions are computed.
{Otherwise}, if \Change{$h/M$ is greater than $\epsilon/|\aaa |$}{$Pe_{h/M}\geq1$}, the algorithm is called recursively, this time requesting the computation of the interactions over elements of size $h/M$.
This recursive algorithm is executed $O(\hspace{0.05 cm} \log (\Change{\hspace{0.05 cm} |\aaa| h \hspace{0.05 cm}/\hspace{0.05 cm}\epsilon }{\Pe_h}))$ times.

\begin{figure}[h!tbp]
    {\scriptsize
    
    \hfil
    \includegraphics[width=.49\textwidth]{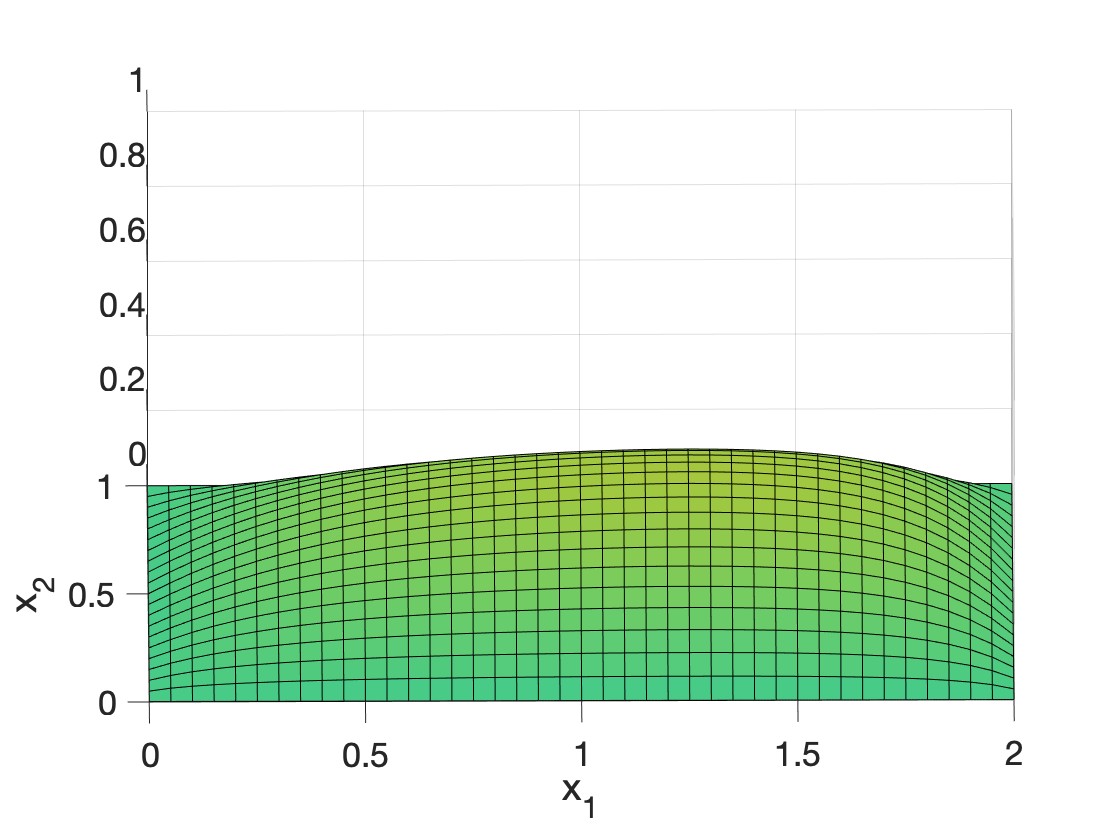}
    \hfil
    \hfil
    \includegraphics[width=.49\textwidth]{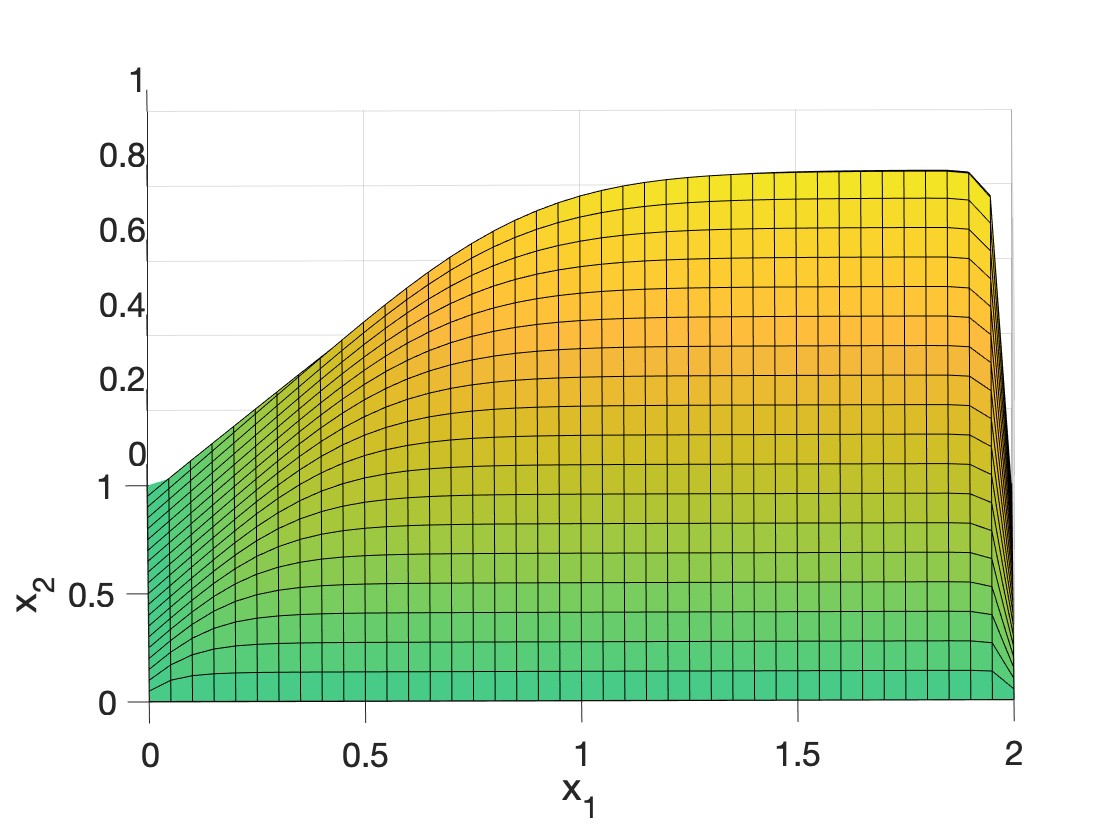}
    \hfil

    \vspace{-2pt}
    \hfil
    $h=0.125\times 10^{-6}$
    \hfil\hfil
    \hfil
    $h=2.5\times 10^{-6}$
    \hfil

    \includegraphics[width=.49\textwidth]{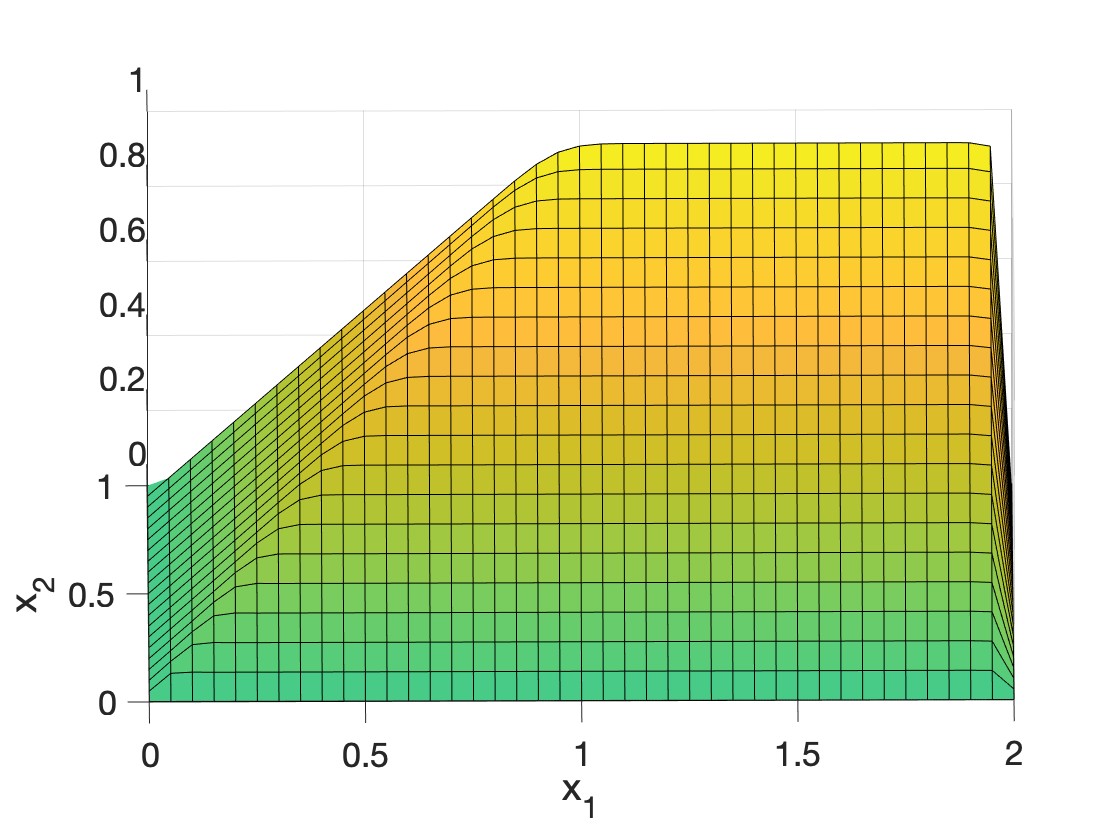}
    \hfil
    \hfil
    \includegraphics[width=.49\textwidth]{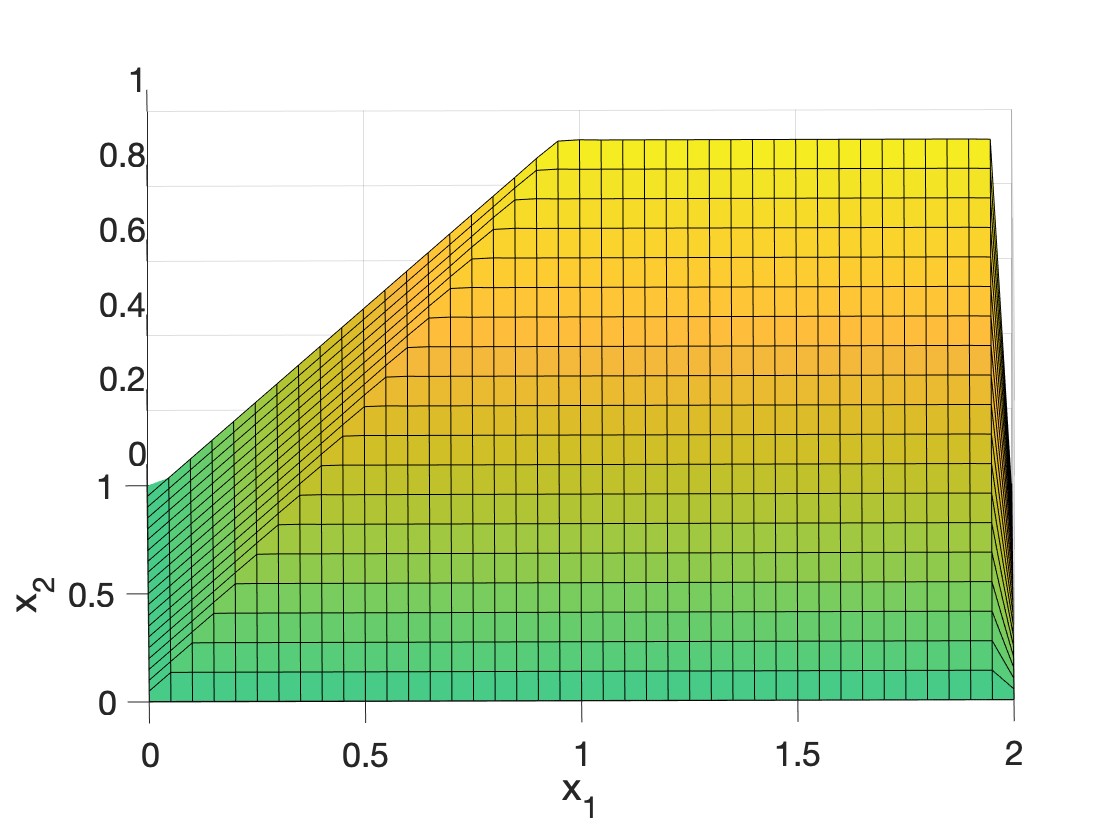}
    \hfil

    \vspace{-2pt}
    \hfil
    $h=0.00005$
    \hfil\hfil
    \hfil
    $h=0.001$
    \hfil

    }
    \caption{\small\label{F:bubbles}  
   Recursive bubbles. Recursive bubbles obtained considering $\aaa=(1,1)$, $c=0$, $\epsilon=10^{-6}$, $f=1$ and $g=0$.
    Smoothest bubble, corresponding to the deepest level in the recursion, where $h<\epsilon$ (top-left), second to last level of recursion (top-right), first level of recursion (bottom-left), actual bubble used in the computation of the solution of the original problem on the physical domain (bottom-right).}
    
\end{figure}

\Change{}{

In the experiments that we show in the rest of this section, we focus in problems with $Pe> 10^3$. 
For the computation of the bubbles, we employ $20 \times 20$ meshes (i.e., $M = 20$). 
The number of recursion steps then depends on the Péclet number of the problem and on the mesh size $h$ of the global discretization. For example, if the Péclet number is $Pe = 10^6$ and the size of the global mesh $h=1/50=0.02$, the finest level is reached in four recursion steps. Indeed, $$\log_{20}(Pe_h) = \log_{20}(h \cdot Pe) = \log_{20} (0.02 \cdot 10^{6} ) \cong 3.31<4;$$ which results in a total of $6 \times 4 = 24$ bubbles being computed overall, namely six bubbles per recursion level: four element bubbles and two patch bubbles. In Figure~\ref{F:bubbles} we show all bubbles corresponding to domains of the type $\omega_S^1$ (two-element patches corresponding to vertical edges) that are computed in the recursive step. It is worth noticing how the bubbles get smoother when the meshsize gets smaller. Finally, it is important to mention that we always rescale the problems and solve on the reference domain $(0,2)\times(0,1)$; computing on the real physical domains proved to be unstable, and rescaling to the reference situation substantially improved our solutions. }

    \delete{
In the experiments that we show in the rest of this section, we use a global mesh with $h=0.02$, $\epsilon = 10^{-6}$ and $1\le |\aaa|\le 2$.  
For the computation of bubbles, we use $20\times20$-meshes (i.e., $M=20$). Therefore, we reach the finest level in four recursion steps, because $\log_{20}h/\epsilon = \log_{20} 0.02/10^{-6} \cong 3.31<4$; thus computing $6\times4=24$ bubbles overall; six bubbles per recursion level: four element- and two patch-bubbles. 
In Figure~\ref{F:bubbles} we show all the bubbles corresponding to domains of the type $\omega_S^1$ (two-element patches corresponding to vertical edges) that are computed in the recursive step. It is worth noticing how the bubbles get smoother when the meshsize gets smaller. Finally, it is important to mention that we always rescale the problems and solve on the reference domain $(0,2)\times(0,1)$; computing on the real physical domains proved to be unstable, and rescaling to the reference situation substantially improved our solutions.
}

\newpage
\Change{}{\subsubsection{Method comparison}}
\new{
In this section, we present examples with and without a reaction term to compare the solutions obtained using the proposed approach with those produced by the traditional residual-free bubble method, which coincides with SUPG for a particular choice of the stabilization parameter. We focus on comparing results obtained with standard element bubbles and with the newly introduced patch bubbles. In Section~\ref{sec:EOC}, where we compute the experimental orders of convergence, we will also include a comparison with a DG approach.
}

\new{
\textit{Example 1: } 
In this example, we compare the solutions obtained with our newly proposed BMZ method, the classical RFB method (RFB/SUPG), and the SUPG method (SUPG/GSC), with the parameter corresponding to GSC from~\cite{Key2023}. 
The latter comparison is included to illustrate that RFB/SUPG and SUPG/GSC not only coincide at the theoretical level, but also produce comparable numerical solutions. The SUPG/GSC results reported here are taken from \cite[Fig.~5.2, Table~5.2]{Key2023}.
}

\new{
We consider an advection-diffusion (AD) problem and an advection-diffusion-reaction (ADR) equation given by Eq.~\eqref{eq:convdiff-stat-strong} with $\Omega$ the unit square, $\epsilon=1$,  $\mathbf{a}=-10^3(\cos(\frac{\pi}{6}),\ \sin(\frac{\pi}{6}))$,  $f=0$, and Dirichlet boundary data given by
\begin{equation*}
    g(x_1,x_2) =
    \begin{cases}
        1,&  \text{ if } x_1=1 \text{ or } x_2=0,\\
        0,& \text{ otherwise. }
    \end{cases}
\end{equation*}
For the reaction term, we take $c=0$ and $c=7500$, respectively. In these problems, the Péclet number is $\Pe = |\aaa|/\epsilon = 10^3$.
It is worth observing that $g$ is discontinuous and does not belong to $H^{1/2}(\partial\Omega)$, and thus lies beyond the regime in which the solution is known to be in $H^1(\Omega)$. Nevertheless, the results are very satisfactory in the internal layer, which starts at the point $\mathbf{x}=(1,1)$ and is aligned with the direction $\aaa$.
}

\begin{figure}[h!tbp]
        \hfil
    \includegraphics[width=.49\textwidth]{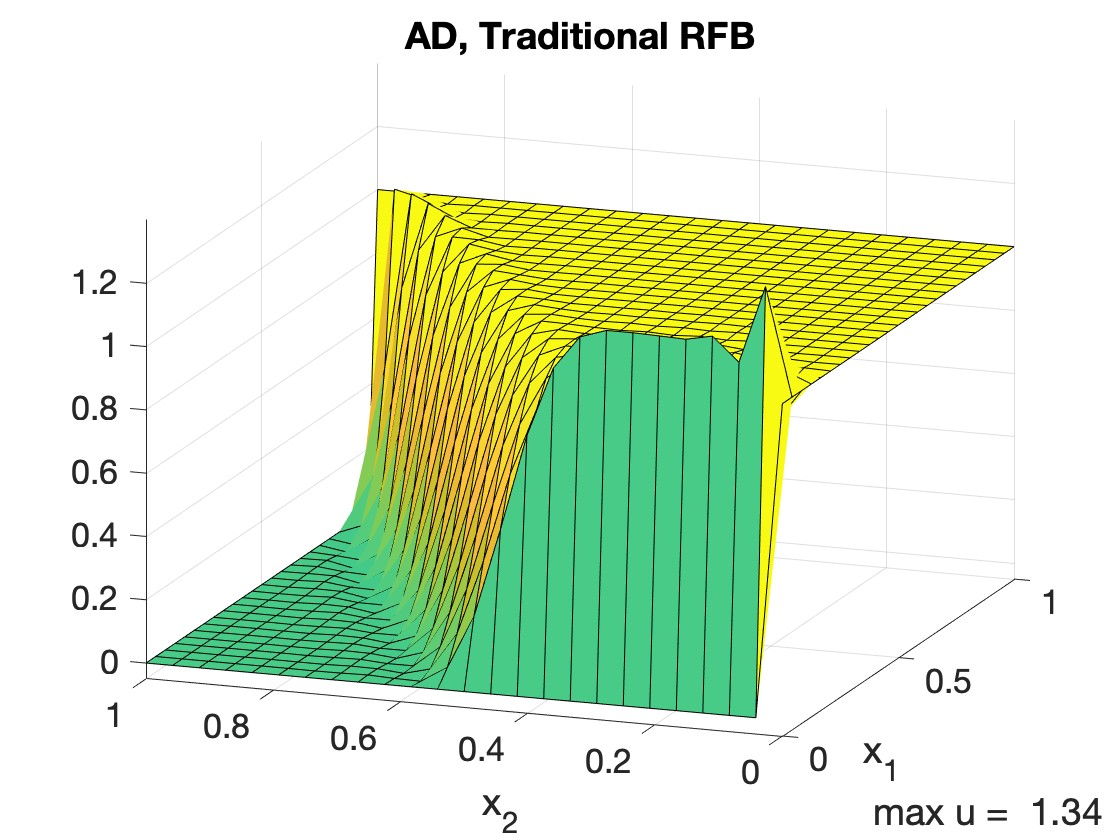}
    \hfil
    \includegraphics[width=.49\textwidth]{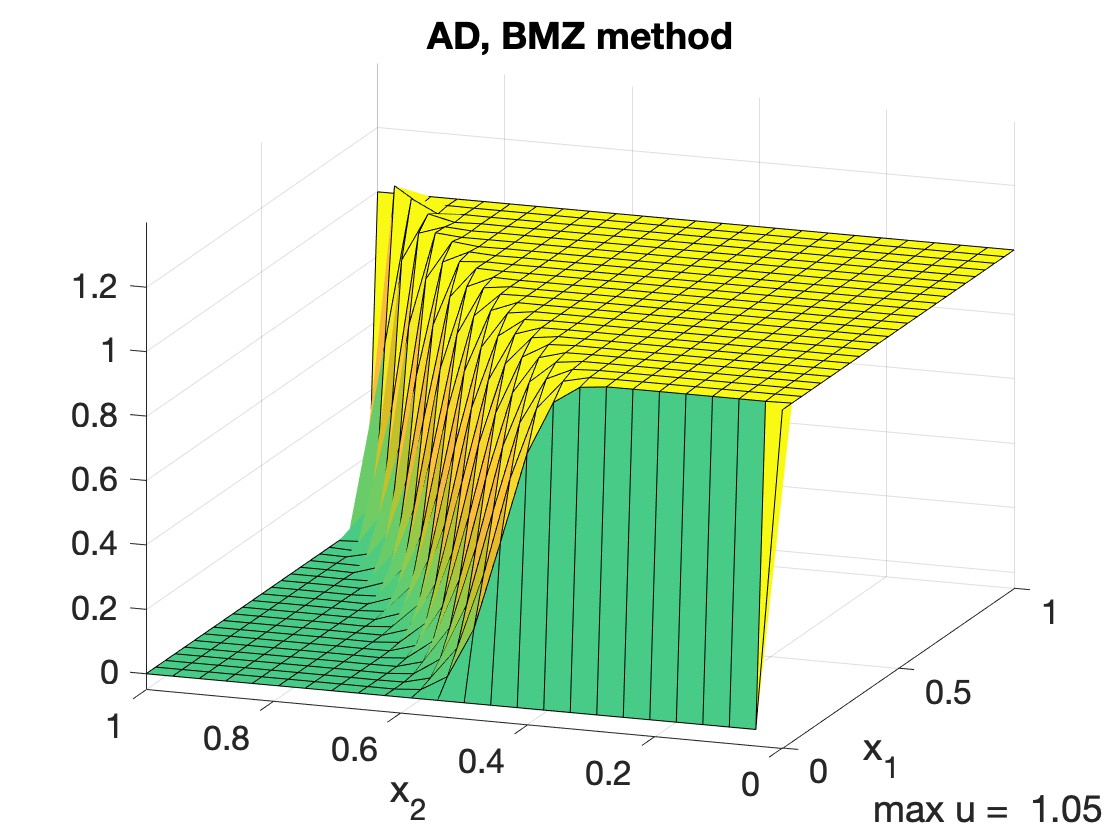}
    \hfil
   \includegraphics[width=.49\textwidth]{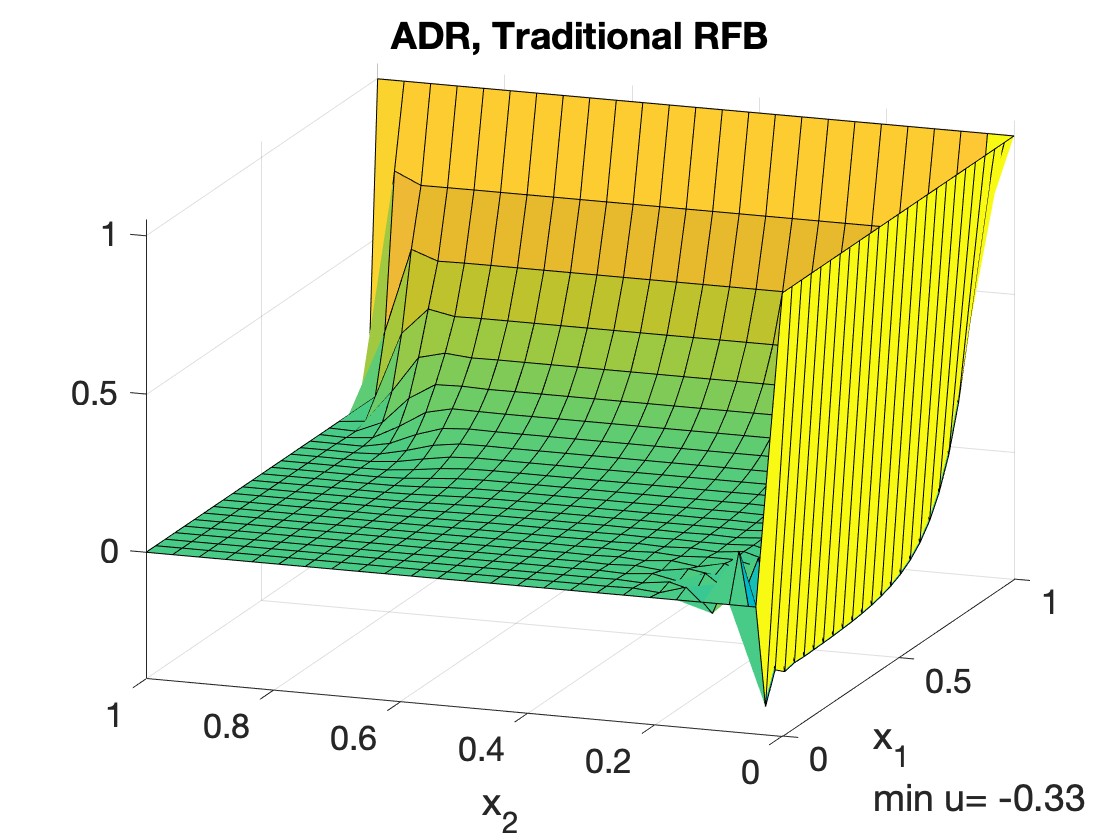}
    \hfil
    \includegraphics[width=.49\textwidth]{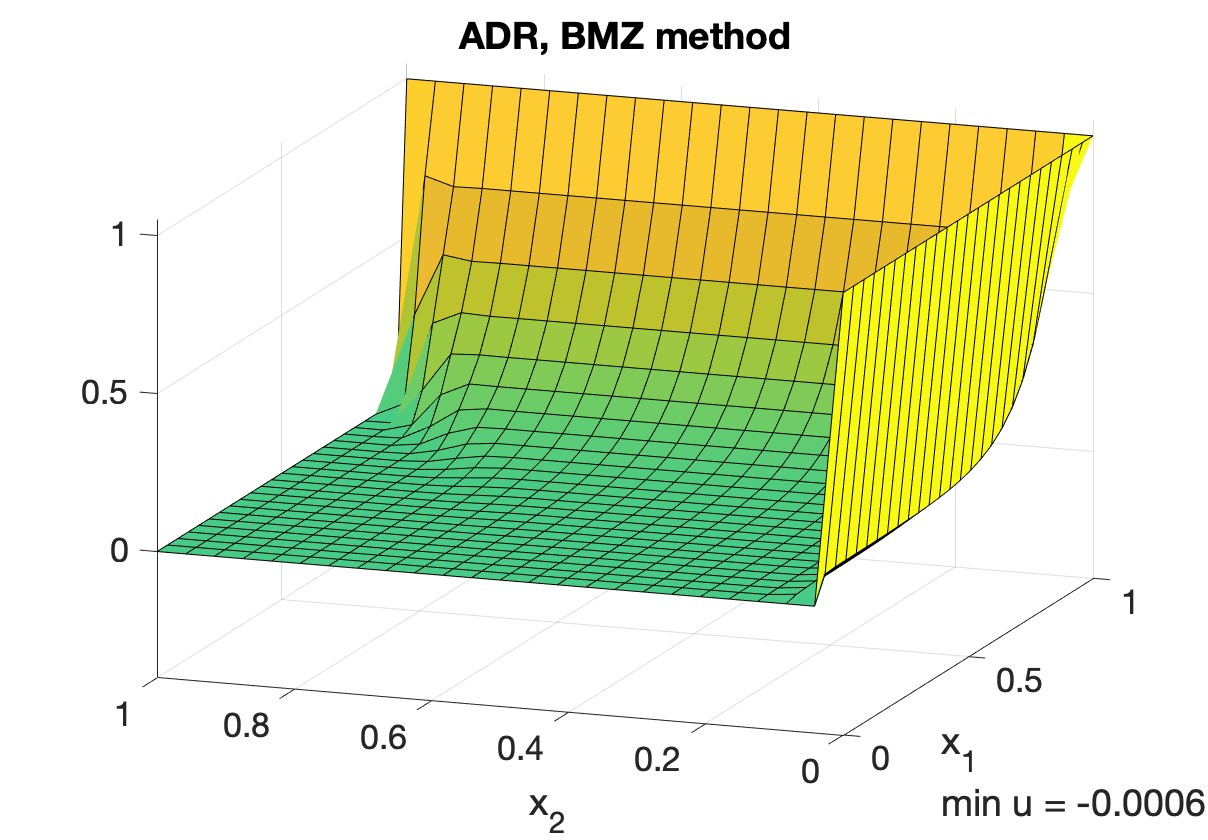}

    \caption{\small\new{
    \label{F:ex-stat-r2} \textit{Example 1,} Comparison of stabilization using the usual residual-free bubbles (left) and BMZ method (right). For the advection-diffusion problem (above) and advection-diffusion-reaction problem (below). The solutions correspond to Eq.~\eqref{eq:convdiff-stat-strong} with $Pe=10^{3}$ with $\epsilon = 10^{-6}$, $\aaa=-10^3(\cos(\frac{\pi}{6}), \sin(\frac{\pi}{6}))$, $f=0$, boundary conditions $u=1$ on $x_1=1$ and $x_2=0$, and a homogeneous Dirichlet condition elsewhere in a partition  with $N\times N$ elements, with $N= 24$. In the ADR problem, the reaction term is $c=7500$. The solution derived using the proposed approach effectively eliminates the oscillations near the boundary layer.}}
\end{figure}

\new{
Figure \ref{F:ex-stat-r2} presents the solutions obtained for RFB/SUPG and BMZ with a partition with $24\times 24$ square elements ($\Pe_h \cong 42$). The solution produced by BMZ effectively eliminates the oscillations near the boundary layer observed in the RFB/SUPG solution, for both the AD and ADR problems. Moreover, the solution obtained with the BMZ method exhibits a satisfactory global behavior, reinforcing its reliability for this problem. The solution with SUPG/GSC can be observed in~\cite[Fig.~5.2]{Key2023}.
}

\new{
    Table~\ref{table: max_min r-2} reports the maximum and minimum values obtained for both problems with each method. The similarity between the values in the second and third rows reinforces the observation that the RFB/SUPG method coincides with the SUPG/GSC method, with a particular choice of parameters; the latter seems to perform better in this case.
For the AD problem with $c=0$, the BMZ method shows an improvement of approximately one order of magnitude in both the minimum and maximum values when compared with the other two methods.
In the ADR case, none of the three methods exceed the prescribed maximum boundary value. However, the error of the minimum value obtained with the BMZ method is $6.37\times 10^{-4}$, which is significantly smaller (three orders of magnitude) than the values $0.3352$ and $0.2805$ obtained with the RFB/SUPG and SUPG/GSC methods, respectively.
}

\begin{table}[ht]
\centering
\new{
\begin{tabular}{|c|c c|c c|}
\hline
Method & 
AD &  & ADR &  \\
 & min & max & min& max  \\
\hline
BMZ  & -0.0037 & 1.047  & -0.0006 & 1 \\
RFB/SUPG  & -0.0159 & 1.342  & -0.3352 & 1 \\
SUPG/GSC & -0.0159 & 1.126  & -0.2805 & 1 \\
\hline
\end{tabular}
}
\caption{\small
\new{
\label{table: max_min r-2} 
\textit{Example 1,} Maximum and minimum values obtained for the AD and ADR problems using the BMZ, RFB/SUPG, and SUPG/GSC methods. The BMZ method yields the most accurate results, both with and without the reaction term, for the maximum and minimum values.}
}
\end{table}

\subsubsection{Experimental order of convergence\label{sec:EOC}}

In the following example, the experimental order of convergence is analyzed. A sequence of uniform meshes with meshsize $h=1/N$, 
is considered, where $N = 10, 20, 40, 80, 160$. 
We present first the results for the
standard $L^2$ and $L^1$ norms, and the three methods: RFB/SUPG, DG and BMZ.

The computations were performed using our own code for RFB/SUPG and BMZ, whereas the toolbox \textsc{festung} {\cite{FESTUNG}} was used for the DG computations, using linear discontinuous finite elements with interior penalization (IIPG), parameter $\eta=6/h$, Zalesak Flux limiters with parameters $\gamma=0$, $\beta=10^{-4}$. 
This choice led, in our experience, to the best results.

In the case of the methods involving bubbles, for the computation of the quadrature in the integrals associated with each norm, all degrees of freedom are taken into account. 
Specifically, according to~\eqref{def:u_h}--\eqref{def:u_h=uL+B+S}, the analyzed difference is given by:  
\new{
\begin{equation*}
    u - u_h = u - (u_L + u_B + u_S)
    \quad\text{and}\quad
    u - u_h = u - (u_L + u_B),
\end{equation*}    
for BMZ and RFB, respectively.
}
Traditionally, only the term $u - u_L$ is considered, which evaluates the difference based on values at the nodes of the partition of the domain $\Omega$. However, incorporating $ u_B $ and $ u_S $ allows for a more accurate representation, as it includes  data of the bubbles from within the elements. This refinement leads to improved estimations in the computation of the norms.
\Change{}{It is important to note that the number of degrees of freedom depends on the chosen method. For a partition consisting of $N \times N$ square elements, the DG method has $4N^2$ degrees of freedom. The RFB and BMZ methods additionally include degrees of freedom associated with the bubble functions. As a result, the RFB method has $5N^2 + 2N + 1$ degrees of freedom, while the BMZ method has $7N^2 - 2N + 3$ degrees of freedom.}
\Change{}{Table~\ref{table:DOFs} shows the number totally available of degrees of freedom associated to each method, depending of the size of the partition $N$.
}

\begin{table}[ht]
\centering
\begin{tabular}{|c|c|c|c|}
\hline
$N$ & BMZ & RFB/SUPG & DG  \\
\hline
10   &    683 &     521 &     400\\
20   &   2763 &    2041 &    1600 \\
40   &  11123 &    8081 &    6400 \\
80   &  44643 &   32161 &   25600 \\
160  & 178883 &  128321 &  102400 \\
\hline
\end{tabular}
\caption{Degrees of freedom corresponding to BMZ, RFB and DG method, depending on the size of the partition $N$.
\label{table:DOFs}
}
\end{table}

\vspace{0.5 cm}
\textit{Example 2: } 
Let us consider the problem \eqref{eq:convdiff-stat-strong} defined on the unit square with $\epsilon = 10^{-6}$, $\aaa = (1,1)$, $c=0$ and homogeneous Dirichlet boundary conditions.
\delete{
\begin{equation*}
    \begin{cases}
        -\epsilon \Delta u + u_x + u_y = f& \text{in $\Omega$},\\
        \hspace{2.3 cm}u = 0 & \text{on $\partial \Omega$}, 
    \end{cases}
\end{equation*}}
\Change{}{Then, the Péclet number $Pe\approx10^6$.} The function on the right-hand side $f$ is defined in such a way that the exact solution is given by
\begin{equation*}
    u(x,y) = 2(\sin x)(1-e^{-(1-x)/\epsilon})y^2(1-e^{-(1-y)/\epsilon}). 
\end{equation*}

\begin{table}[ht]
\centering
\begin{tabular}{|c|c|c|c|c|c|c|}
\hline
$N$ & BMZ & EOC & RFB/SUPG & EOC & DG & EOC \\
\hline
10  & 2.250\,e\,-3  &  -  & 1.176\,e\,-1 &  -  & 0.854\,e\,-3 &  -  \\
20  & 0.564\,e\,-3  & 1.9 & 0.850\,e\,-1 & 0.4 & 0.221\,e\,-3 & 1.9  \\
40  & 0.143\,e\,-3  & 1.9 & 0.608\,e\,-1 & 0.4 & 0.057\,e\,-3 & 1.9  \\
80  & 0.048\,e\,-3  & 1.5 & 0.432\,e\,-1 & 0.4 & 0.020\,e\,-3 & 1.4  \\
160 & 0.045\,e\,-3  & 0.0 & 0.306\,e\,-1 & 0.4 & 0.020\,e\,-3 & -0.0 \\
\hline
\end{tabular}
\caption{\small
    $L^2(\Omega)$ errors and experimental order of convergence. 
    We consider a homogeneous Dirichlet problem, with $\epsilon =10^{-6}$.
    It is worth observing that the BMZ and the DG method exhibit a second order reduction, at least in the pre-asymptotic regime,
    whereas the other one (RFB/SUPG) seems to converge with order $0.4$.
    It is also striking that the errors are already of order $10^{-3}$ for $N=10$ for the BMZ and DG methods, and also at least two orders of magnitude smaller than the other one, for some values of $N$. The errors for the DG method are smaller than those for BMZ by a factor $1/2$, at the expense of having discontinuous elements though.}
\label{table: EOC L2}
\end{table}

\begin{table}[ht]
\centering
\begin{tabular}{|c|c|c|c|c|c|c|}
\hline
$N$ & BMZ & EOC & RFB/SUPG & EOC & DG & EOC \\
\hline
10  & 1.810\,e\,-3 &  -   & 3.776\,e\,-2 &   -  & 0.689\,e\,-3 &   -  \\
20  & 0.453\,e\,-3 & 1.9  & 1.957\,e\,-2 & 0.9  & 0.181\,e\,-3 & 1.9  \\
40  & 0.114\,e\,-3 & 1.9  & 1.013\,e\,-2 & 0.9  & 0.046\,e\,-3 & 1.9  \\
80  & 0.029\,e\,-3 & 1.9  & 0.517\,e\,-2 & 0.9  & 0.012\,e\,-3 & 1.8  \\
160 & 0.008\,e\,-3 & 1.7  & 0.261\,e\,-2 & 0.9  & 0.004\,e\,-3 & 1.5  \\
\hline
\end{tabular}
\caption{\small
    $L^1(\Omega)$ errors and experimental order of convergence. 
    We consider a homogeneous Dirichlet problem, with $\epsilon =10^{-6}$.
    It is worth observing that the BMZ and the DG method exhibit a second order reduction, at least in the pre-asymptotic regime,
    whereas the other one (RFB/SUPG) seems to converge with order $0.9+$.
    It is also striking that the errors are already of order $10^{-3}$ for $N=10$ for the BMZ and DG methods, and also at least two orders of magnitude smaller than the other one, for some values of $N$. The errors for the DG method are smaller than those for BMZ by a factor $1/2$, at the expense of having discontinuous elements though.}
\label{table: EOC L1}
\end{table}

    In Table~{\ref{table: EOC L2}} (Table~{\ref{table: EOC L1}}) we show the $L^2(\Omega)$ ($L^1(\Omega)$) norm of the errors corresponding to the three methods. 
It is worth observing that the BMZ and the DG methods exhibit a second order reduction of the $L^2$ and $L^1$ norms of the error, whereas RFB/SUPG exhibits an order around $0.4$ and $0.9$, respectively.
It is important to note, however, that to the best of our knowledge, no theoretical error estimates have been developed for the $L^2$ and $L^1$ norms for any of the methods considered in this article.
It is striking that the errors in the $L^1$ and $L^2$ norms are already of order $10^{-3}$ for $N=10$ for the BMZ and the DG method, and also at least two orders of magnitude smaller than RFB, for some values of $N$.
We observe that the DG method performs better than BMZ by a factor 2, when measuring the error in the $L^2$ and $L^1$ norms for $\epsilon = 10^{-6}$. 
This good performance of DG comes at the expense of a poor robustness, though, and we discuss this at the end of the section.
\new{We finally notice that the errors for BMZ and DG do not decrease when going from $N=80$ to $N=160$. A possible explanation for this behavior of BMZ stems from considering the theoretical results from Cangiani-Süli~\cite{Cangiani2005long}, which suggest that the error for BMZ would be of order $\sqrt{\epsilon/h}+h$. This bound decreases for certain range of values of $h$ but does not go to zero when $h\to 0$. We believe that we have reached a parameter region where $\text{Pe}_h\not\gg 1$ and the error stays around 1e-5 for a while, until $h$ reaches the region where $\text{Pe}_h < 1$, where the Galerkin solution would converge to zero as expected. Something similar might be happening for DG.}

\begin{table}[ht]
\centering
\begin{tabular}{|c|c|c|c|c|c|c|}
\hline
$N$ & BMZ & EOC & RFB/SUPG & EOC & DG & EOC \\
\hline
10  & 8.078\,e\,-2 &  -   & 1.364\,\,&  -   & 4.634\,e\,-2 & - \\
20  & 5.467\,e\,-2 & 0.5  & 1.487\,\,& -0.1 & 3.304\,e\,-2 & 0.4\\
40  & 2.960\,e\,-2 & 0.9  & 0.821\,\,&  0.8 & 1.925\,e\,-2 & 0.7\\
80  & 1.604\,e\,-2 & 0.9  & 0.874\,\,& -0.0 & 1.034\,e\,-2 & 0.8\\
160 & 0.841\,e\,-2 & 0.9  & 3.084\,\,& -1.8 & 0.535\,e\,-2 & 0.9\\
320 & 0.424\,e\,-2 & 1.0  & 4.260\,\,& -0.4 & - & - \\
640 & 0.211\,e\,-2 & 1.0  & 1.334\,\,&  1.6 & - & - \\
\hline
\end{tabular}
\caption{\small
    $H^1$ error in $(2h,1-2h)^2$.
We consider a homogeneous Dirichlet problem, with $\epsilon =10^{-6}$.
    In order to avoid the boundary layer, we compute the error only on the subdomain $[2h,1-2h]^2$.
    The BMZ and DG methods behave very similarly, at least for $N \le 160$, with errors two orders of magnitude smaller than those of RFB/SUPG; \textsc{festung} gives a warning and halts for $N\ge 200$.
    It is worth observing that the BMZ method exhibits a steady decrease of the error of order $1$, and the same seems to happen for DG, whereas the RFB method behaves more erratically.
    }
\label{table: EOC H1}
\end{table}

\begin{table}[ht]
\centering
\begin{tabular}{|c|c|c|c|c|}
\hline
$N$ & BMZ & EOC & RFB/SUPG & EOC \\
\hline
10  & 2.1481\,e\,-2 &  -   & 0.3628 &  -   \\
20  & 1.0281\,e\,-2 &  1.0 & 0.2797 &  0.3 \\
40  & 0.3936\,e\,-2 &  1.3 & 0.1092 &  1.3 \\
80  & 0.1508\,e\,-2 &  1.3 & 0.0821 &  0.4 \\
160 & 0.0559\,e\,-2 &  1.4 & 0.2050 & -1.3 \\
320 & 0.0199\,e\,-2 &  1.4 & 0.2003 &  0.0 \\
640 & 0.0070\,e\,-2 &  1.5 & 0.0443 &  2.1 \\
\hline
\end{tabular}
\caption{\small
    Stability norm in $(2h,1-2h)^2$.
We consider a homogeneous Dirichlet problem, with $\epsilon =10^{-6}$, only for the two methods which work with continuous elements.
    In order to avoid the boundary layer, we compute the error only on the subdomain $[2h,1-2h]^2$.
    It is  worth observing that the BMZ method exhibits a steady decrease of the error of order $1.5$, whereas the RFB/SUPG method behaves more erratically.
    }
\label{table: EOC stab-norm} 
\end{table}

In Table~{\ref{table: EOC H1}} we show the errors measured in the (broken) $H^1$-seminorm. It is worth observing that BMZ and DG outperform RFB/SUPG by two orders of magnitude. Also, notice that for DG the approximations do not belong to $H^1$ and we compute these solutions up to $N=160$; \textsc{festung} gives a warning and halts for $N\ge 200$. It is worth observing that the BMZ method exhibits a steady decrease of the error of order $1$, and the same seems to happen for DG, whereas the RFB method behaves more erratically.

In Table~{\ref{table: EOC stab-norm}} we show the errors for BMZ and RFB/SUPG measured in the following part of the stability norm:

\begin{equation}\label{stab-norm}
    \left(
        \sum_T h_T||\aaa \cdot \nabla  v||_{L^2(T)}^2\right)^{1/2}.
\end{equation}
In {\cite{Brezzi1999}} they defined this norm, in which they obtained error estimates for the RFB/SUPG method.
For this reason, we include this last table to compare the behavior of BMZ and RFB/SUPG with respect to this norm.

In order to avoid the boundary layers we have computed these norms involving gradients of the solution on the subdomain $[2h,1-2h]^2$.
\delete{We do not present here the $H^1$ seminorm and the stability norm of the error on the whole domain, because the quadrature on the considered meshes is not sufficient to measure them efficiently.
To illustrate this difficulty, it is worth observing that the $H^1(\Omega)$ seminorm of the exact solution $u$, computed with such quadratures, for $N=10$, $20$, $40$, $80$, $160$, $320$ is approximately
$45221.63$,
$32079.54$,
$22720.10$,
$16078.43$,
$11373.72$,
$8044.04$,
respectively.
This shows that we are still far from a precise computation, which is necessary to observe the error reduction and the experimental order of convergence.}
It is  worth observing that the BMZ method exhibits a steady decrease of the error of order $1.5$, whereas the RFB/SUPG method behaves more erratically.

Summarizing, the presented tables demonstrate a better behavior of BMZ and DG across all norms when compared to those obtained with the traditional 
RFB method, which coincides with the famous SUPG. The orders of BMZ and DG coincide and in some cases DG shows errors a bit smaller than BMZ.
We also recall that the DG method performs better than BMZ by a factor 2, when measuring the error in the $L^2$ and $L^1$ norms for $\epsilon = 10^{-6}$. 
This good performance of DG comes at the expense of a poor robustness, though. Indeed DG seems to perform not so well for larger values of $\epsilon$, and some tuning appears to be necessary. In order to see this, besides looking at Figure~{\ref{F:DG}}, one can see Table~\ref{table:DG larger epsilon}, where the errors are computed for $\epsilon = 10^{-3}$.



\begin{table}[ht]
    \begin{center}
{
    
    
    \begin{tabular}{|c|c|c|c|c|c|c|c|c|}
    \hline 
    \multicolumn{1}{|c|}{$ $}&
    \multicolumn{1}{c}{$ $} & \multicolumn{1}{r}{$L^2(\Omega)$} & \multicolumn{1}{l}{error}& \multicolumn{1}{c|}{} & 
    \multicolumn{1}{c}{$ $} & \multicolumn{1}{r}{$L^1(\Omega)$} & \multicolumn{1}{l}{error}& \multicolumn{1}{c|}{} \\
    \hline 
    N & BMZ & EOC &  DG & EOC& BMZ & EOC &  DG & EOC\\
    \hline
    10  & 2.91\,e\,-3 &  -   &  05.2\,e\,-3  &   -  & 1.93\,e\,-3 &  -   &  1.93\,e\,-3  &   -  \\
    20  & 2.16\,e\,-3 &  0.4 &  07.0\,e\,-3  & -0.4 & 0.74\,e\,-3 & 1.3  &  1.59\,e\,-3  &  0.2 \\
    40  & 1.76\,e\,-3 &  0.2 &  08.3\,e\,-3  & -0.2 & 0.32\,e\,-3 & 1.1  &  1.31\,e\,-3  &  0.2 \\ 
    80  & 1.06\,e\,-3 &  0.7 &  07.1\,e\,-3  &  0.2 & 0.12\,e\,-3 & 1.4  &  0.77\,e\,-3  &  0.7 \\ 
    160 & 0.41\,e\,-3 &  1.3 &  10.1\,e\,-3  & -0.5 & 0.03\,e\,-3 & 1.7  &  0.93\,e\,-3  & -0.2 \\  
    \hline
    \end{tabular}
    }
\end{center}
\caption{\small $L^2(\Omega)$ and $L^1(\Omega)$ error for the BMZ and DG method. We consider $\epsilon=10^{-3}$. The BMZ method achieves smaller errors and higher EOC than the DG scheme in both norms, showing a monotone error decay in contrast with the oscillatory behavior of the DG method.}
\label{table:DG larger epsilon}
\end{table}


\section{Unsteady Problems}
\label{S:RFB-instat}

We now consider the unsteady advection dominated advection-diffusion problem
\begin{equation}\label{eq:inst-problem}
    \left\{
\begin{aligned}
    u_t - \epsilon \Delta u + \aaa\cdot \nabla  u + c\,u&=f \qquad& &\text{in $\Omega$,\ \ $t>0$,}\\
        u &= g \qquad& &\text{on $\partial\Omega$},\ \ t>0,\\
    u(\cdot,0) &= u^0 \qquad& & \text{in }\Omega.\\
\end{aligned}
    \right.
\end{equation}

\subsection{New method}



The proposed approach for the temporal discretization is based on retaining the information from the bubble functions rather than applying static condensation. 
This retention of bubble information greatly simplifies the implementation and separates the idea of static condensation, as a second step, from the idea of the bubbles to be used.
Static condensation could also be done, in order to obtain smaller systems to solve in each timestep, and this will be subject of future research.

The bilinear form associated with problem \eqref{eq:inst-problem} is given by:
\begin{equation*}
a(u,v) = \int_{\Omega} \epsilon\nabla  u\cdot\nabla  v + \aaa \cdot \nabla  u v + c\,uv\hspace{0.1 cm} dx
\end{equation*}

Thus, the weak form of the unsteady problem can be expressed as:
\begin{equation}\label{eq:inst-weakform}
\text{Find } u(t,\cdot) \in \Hoi:\qquad
	(u_t,v) + a(u,v) = (f,v) \qtfa v \in \Hoi.
\end{equation}

After space discretization, the resulting matrix system takes the form:
\[
M U_t + A U = F,
\]
where $M$ represents the mass matrix and $A$ correspond to the matrix associated with the bilinear form $a(\cdot,\cdot)$, incorporating both finite element basis functions and element- and patch-supported bubble functions,
so that $U(t)$ is the vector of all coefficients appearing in~\eqref{def:u_h=uL+B+S} of the discrete unknown $u_h(t)$ at time $t$.

Subsequently, standard time discretization schemes can be applied, such as the Backward Euler, Crank-Nicolson, or Backward Differentiation Formulas, each offering distinct characteristics. 
In the computations below we have chosen to use Crank-Nicolson's scheme, i.e., 
%
%
%
\begin{equation*}
	(M + \frac{\Delta t}2  A ) U^{n} = (M-\frac{\Delta t}2 A) U^{n-1} + \Delta t F_{n-\frac12},
\end{equation*}


The following section demonstrates how this idea leads to a stable and accurate method. In particular, the proposed method effectively eliminates oscillations, achieves the expected maximum and minimum values, and enhances overall performance.

\subsection{Numerical Experiments}

\textit{Example 3: Constant velocity unsteady problem.}
Consider the problem~\eqref{eq:inst-problem}
with $\Omega$ an $L$-domain, with $L:=(0,2)\times(0,2)\setminus(1,2)\times(1,2)$, $\epsilon=10^{-6}$,  $\aaa = (-2,1)$, $c=0$, $f=1$, and $g=0$, with $Pe\approx10^{6}$.

\begin{figure}[h!tbp]

    \includegraphics[width=.32\textwidth]{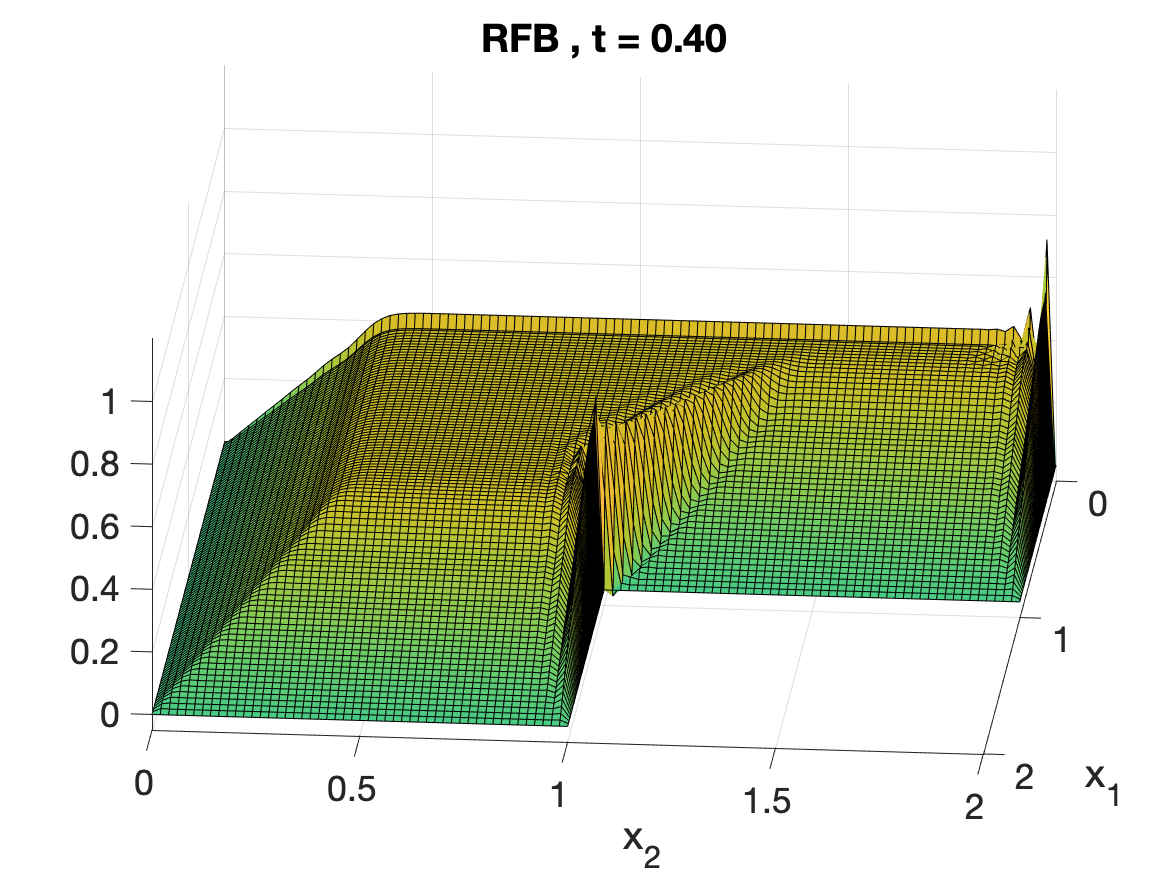
    }
    \includegraphics[width=.32\textwidth]{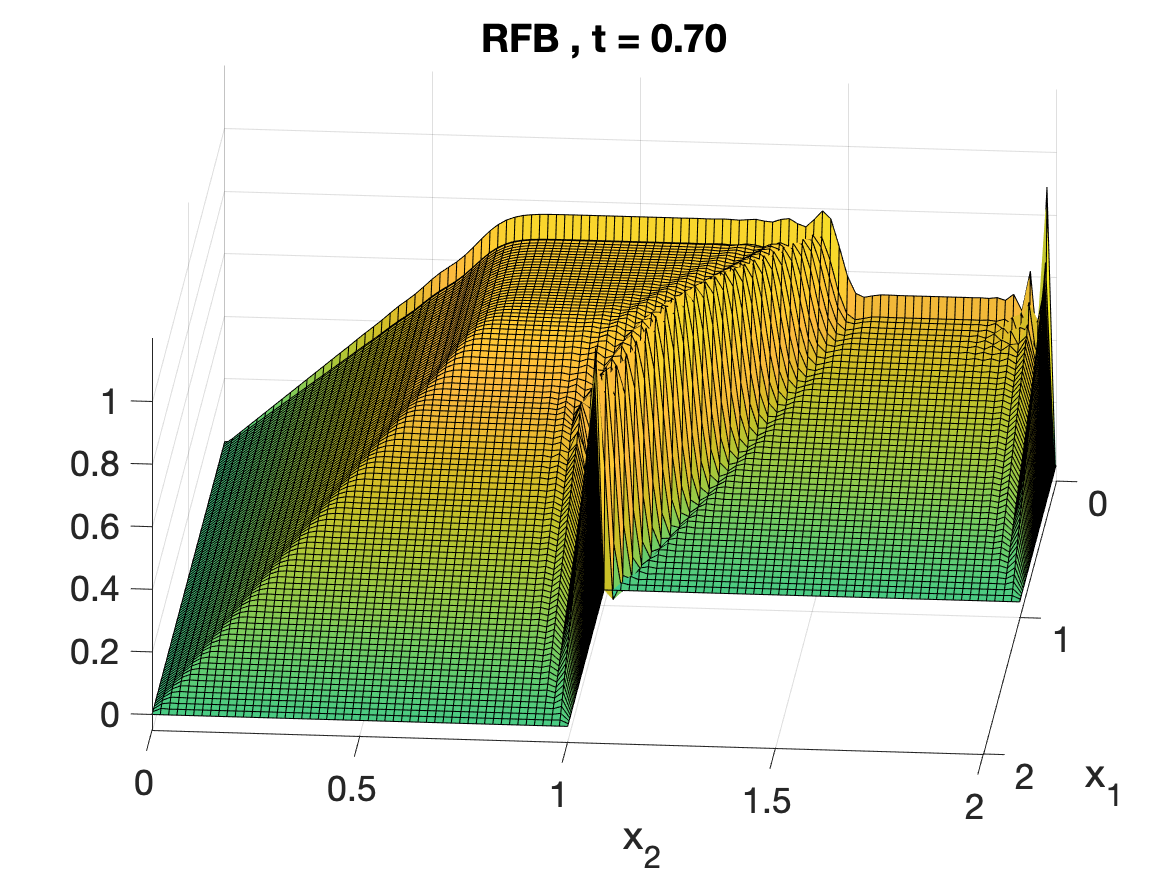}
    \includegraphics[width=.32\textwidth]{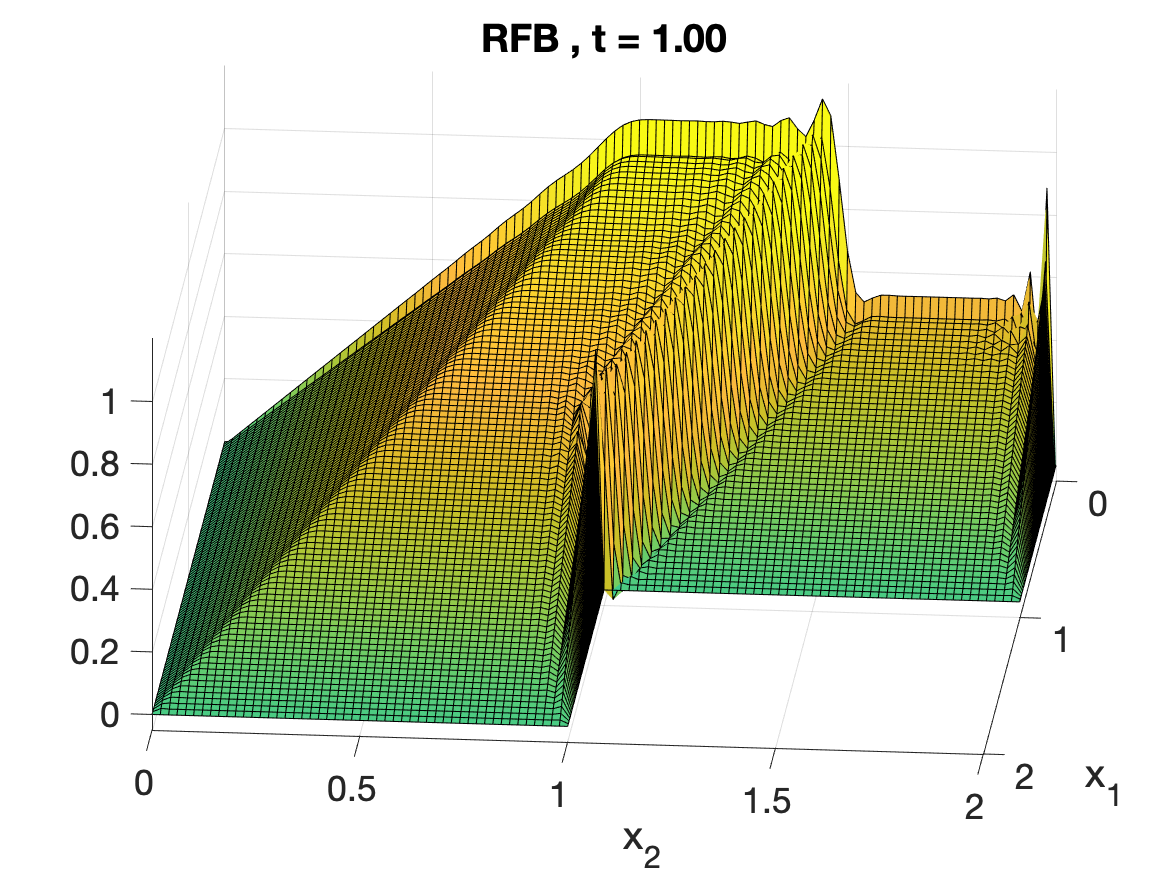}

     \includegraphics[width=.32\textwidth]{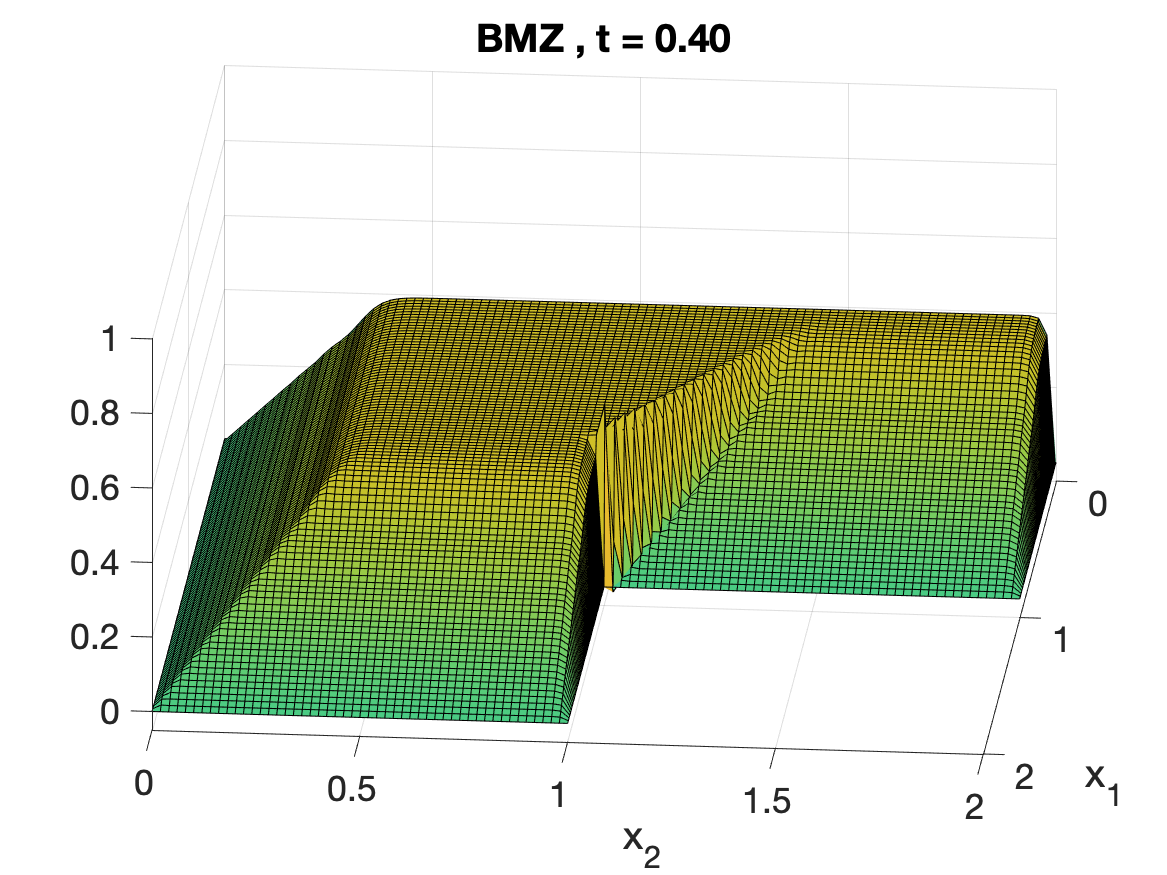
    }
    \includegraphics[width=.32\textwidth]{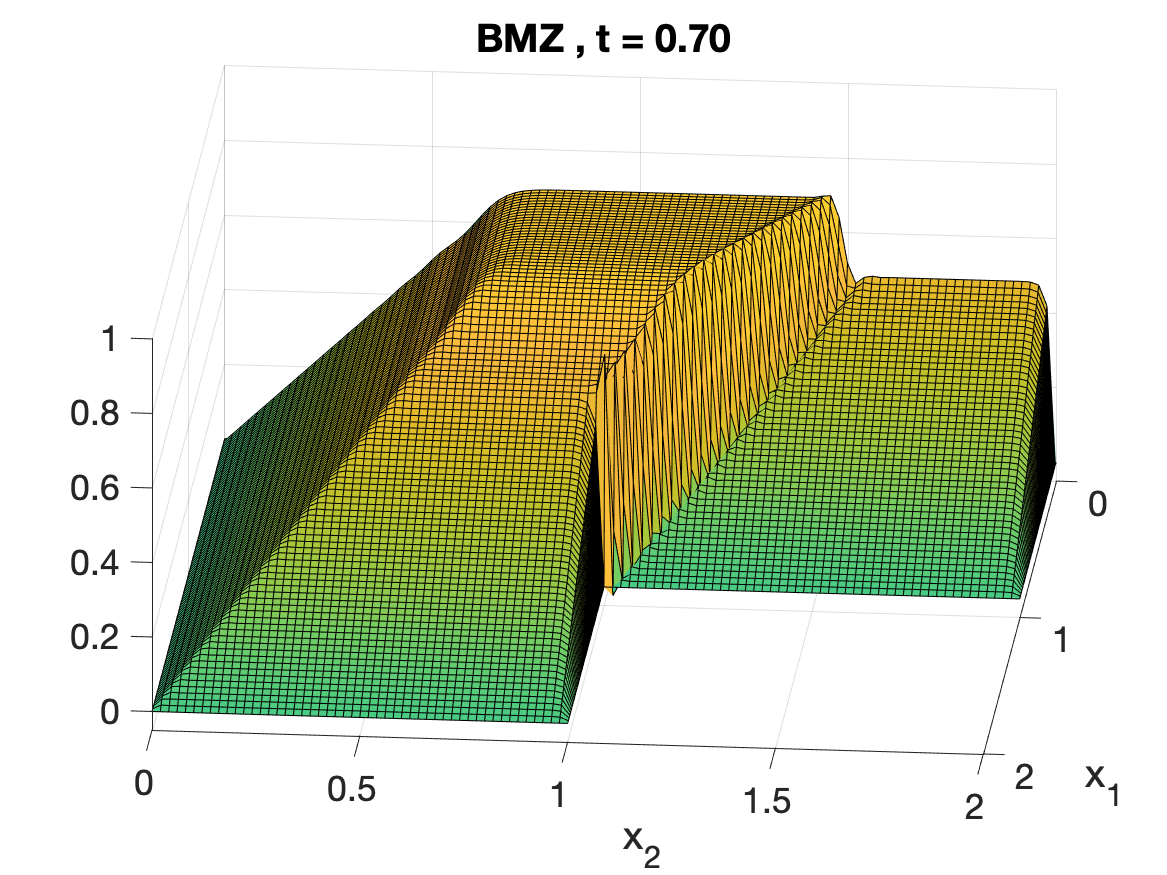}
    \includegraphics[width=.32\textwidth]{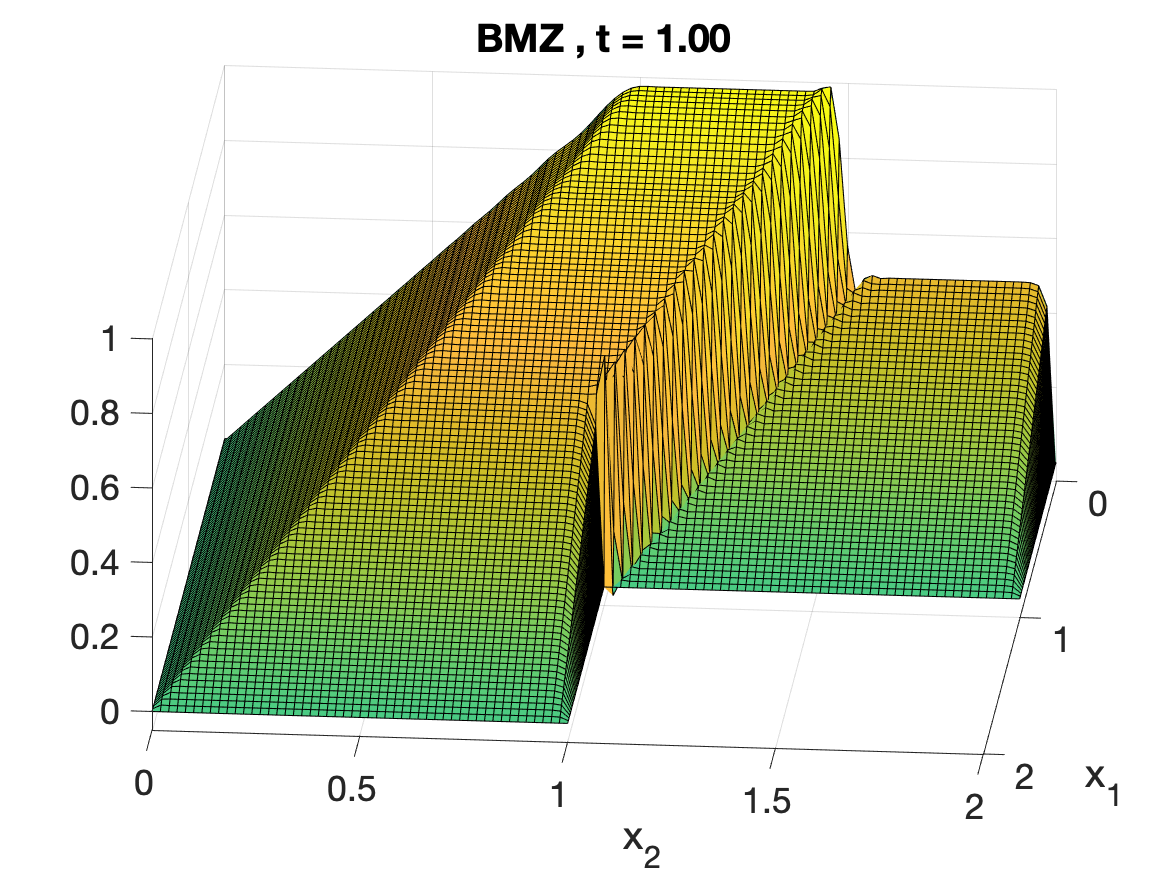}
    
    \caption{\small\label{F:ex4} \textit{Example 3. Constant velocity unsteady problem:} Comparison of stabilization using the traditional residual-free bubbles (above) and BMZ method (below). The solution correspond to \eqref{eq:inst-problem} with $\epsilon=10^{-6}$ for $t= 0.40, 0.70,1$.
    The proposed BMZ method eliminates boundary layer oscillations across all time steps.    
    } 
\end{figure}

As time progresses, the solution grows until reaching $t=1$. Beyond this time, the solution converges to the steady-state problem.
Figure \ref{F:ex4} shows the solutions obtained with the RFB method (above) and the BMZ method (below) at various time steps. The domain is discretized using a mesh with $h=1/50$, and the time step is set as $\Delta t = 0.02$ for the Crank-Nicolson method.

As in the steady case, the solution obtained with the traditional RFB method exhibits oscillations near the boundary and interior layers. This issue is reflected in the maximum values attained at each time. In contrast, the proposed BMZ method delivers accurate results. For instance, at $ t=0.7 $, the expected maximum value is approximately $0.7-\epsilon$.  The BMZ method yields a maximum value of \(0.71722\), while the traditional RFB method reaches \(0.8975\). This discrepancy persists throughout all time steps, highlighting the superior accuracy of the BMZ approach in handling boundary layers.

As a further remark, it is worth highlighting that the implementation of the bubbles defined by Cangiani and Süli (with the recursive computation), combined with our proposed approach for unsteady problems, yields solutions comparable to those obtained with the BMZ method in this example.

\vspace{0.5 cm}
\textit{Example 4: Rotating velocity field.}
We consider problem~\eqref{eq:inst-problem} where
\begin{wrapfigure}[9]{r}{4 cm}
	\vspace{-0.3 cm}
	\includegraphics[width=1\linewidth]{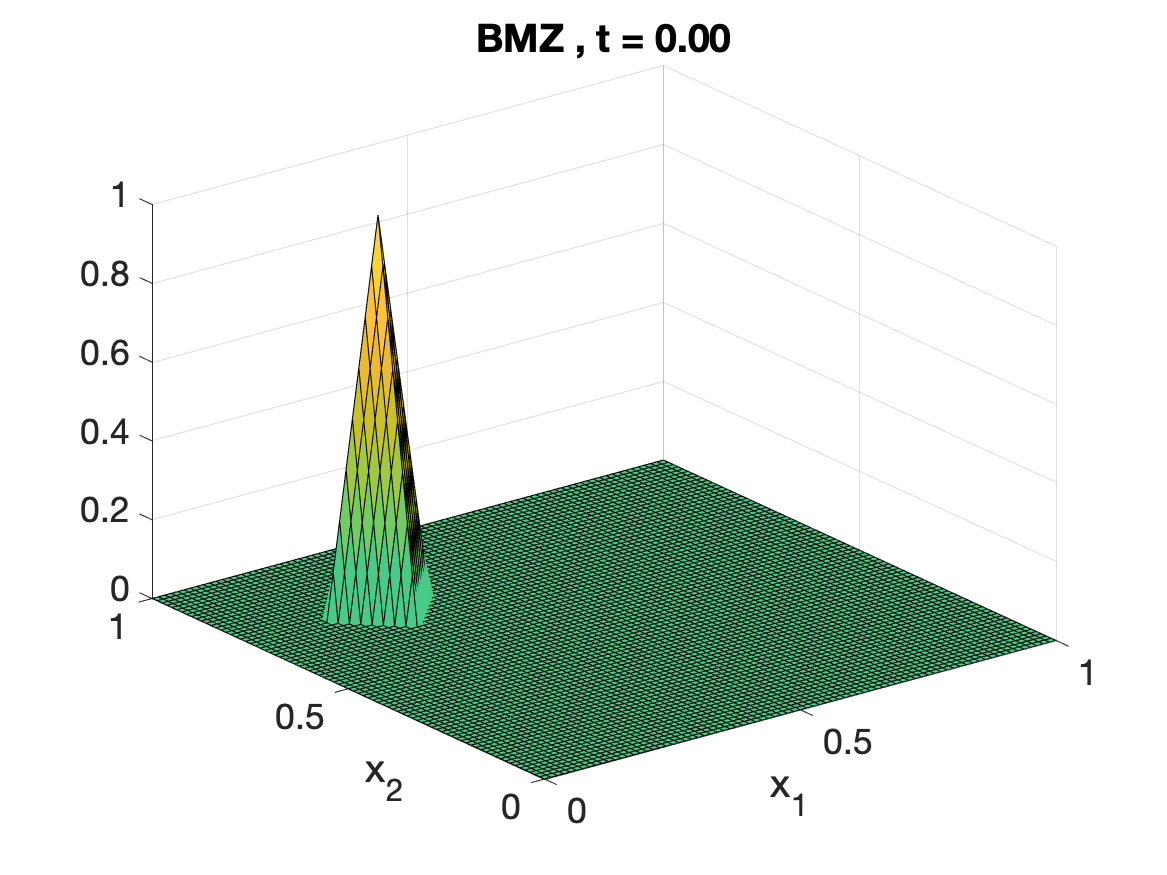}
\caption{\small\label{F:Pyramid}pyramid $u^0 (x).$}
\end{wrapfigure} 
$\Omega$ is the unit square, $\epsilon=10^{-6}$, $c=0$, $f=1$, $g=0$, and \aaa is spatially dependent, representing a clockwise rotation ${\aaa(x_1,x_2) = ( x_2 - 1/2, 1/2 - x_1)}$. The initial function $u^0$, shown in Figure \ref{F:Pyramid}, is defined as the pyramid
  \begin{equation*}
 	u^0(x)= \max\{0,1-\frac{1}{r}||x-x_0||_1\} 
 \end{equation*}
where $x_0 = (0.25 , 0.75) $ and $r = 0.1$.

\begin{figure}[h!]
\centering
{\includegraphics[width=.49\textwidth]{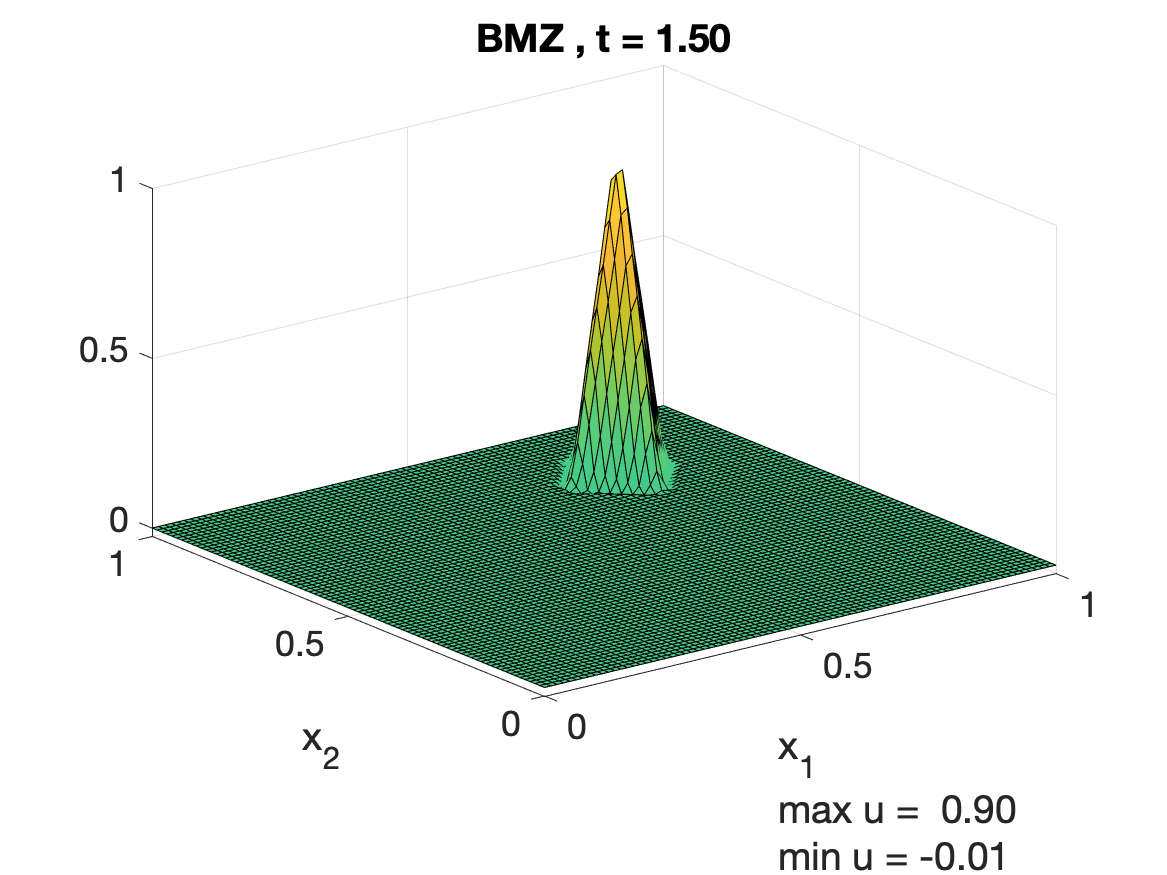}}
\includegraphics[width=.49\textwidth]{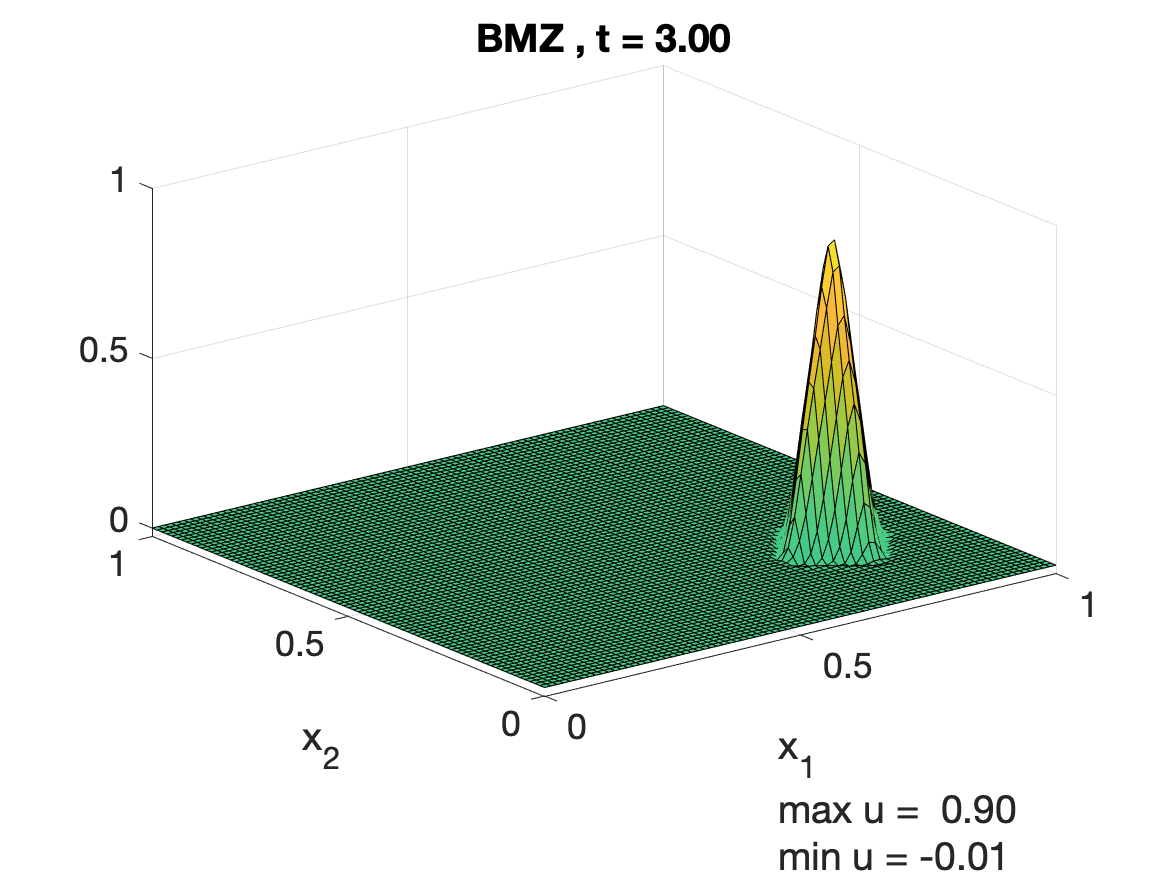}
    \includegraphics[width=.49\textwidth]{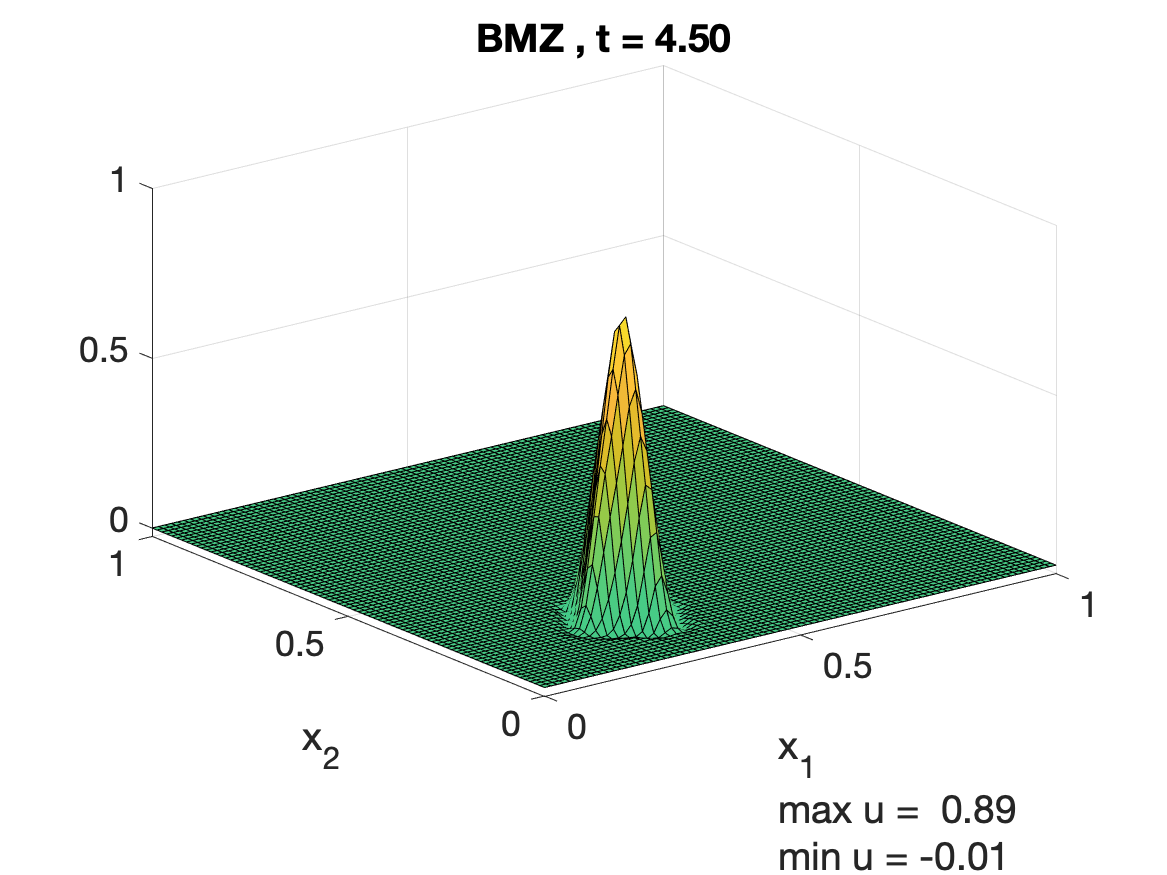}
\includegraphics[width=.49\textwidth]{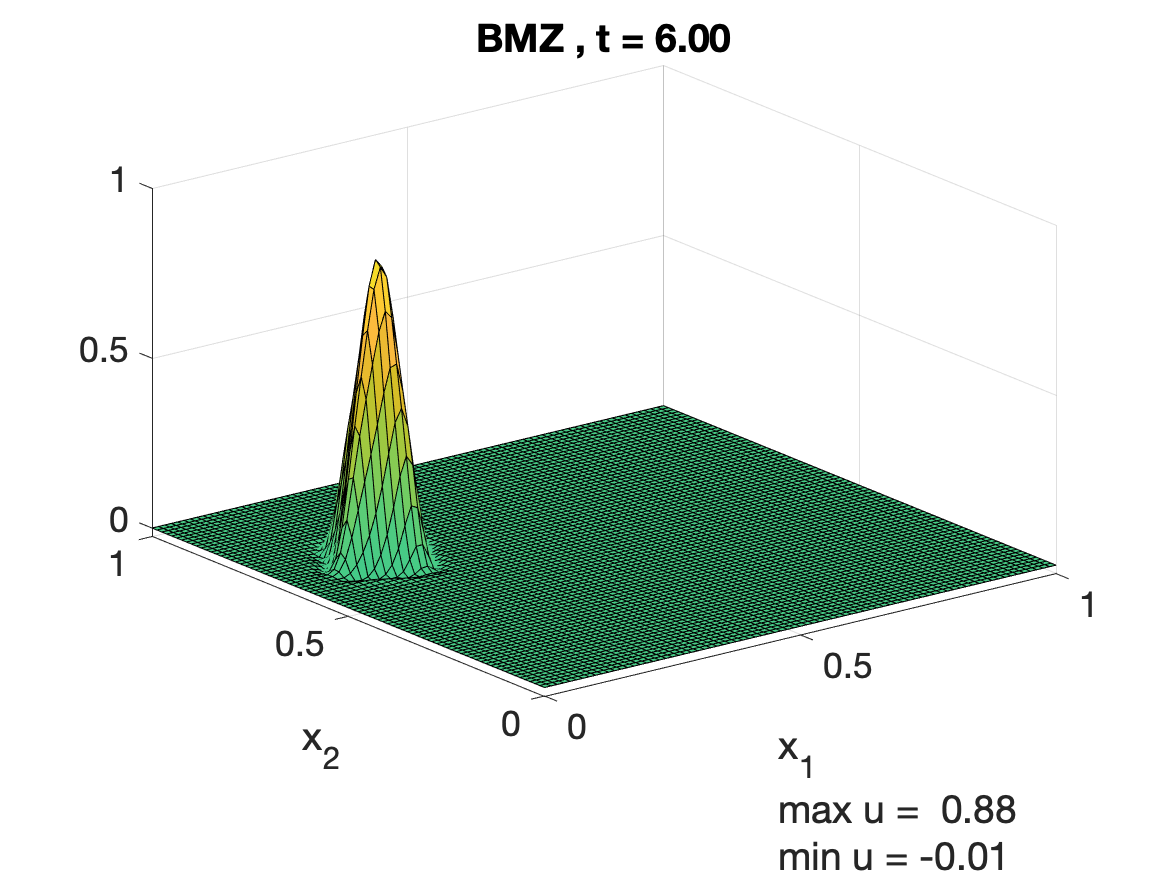}

    \caption{\small\label{F:ex5} \textit{Example 4. Clockwise rotation:} Advection dominated advection diffusion of a pyramid with a clockwise rotation (non constant) velocity.
    The solution corresponds to Eq.~\eqref{eq:inst-problem} with $\epsilon=10^{-6}$ for $t=1.5, 3, 4.5,6$.
    No oscillations are observed at any timestep, and the edges of the pyramid, though slightly diffused, are well-preserved. } 
\end{figure}

Since this example is convection-dominated, because $\Pe\approx10^{6}$, the solution over time closely resembles that of a transport problem. 
Figure \ref{F:ex5} presents the solution obtained using the patch bubbles from the BMZ method. The domain is discretized with a mesh size of $ h = 1/80 $, and the time step is set as $\Delta t = 0.0025 $ for the Crank-Nicolson method.

Our method demonstrates a satisfactory performance across multiple aspects. Firstly, only minor  oscillations are observed. Additionally, the edges of the pyramid, though slightly diffused, are well-preserved. Lastly, the maximum and minimum values of the solutions align closely with the expected ones.

\begin{figure}[h!]
    \centering
    \includegraphics[width=.6\textwidth]{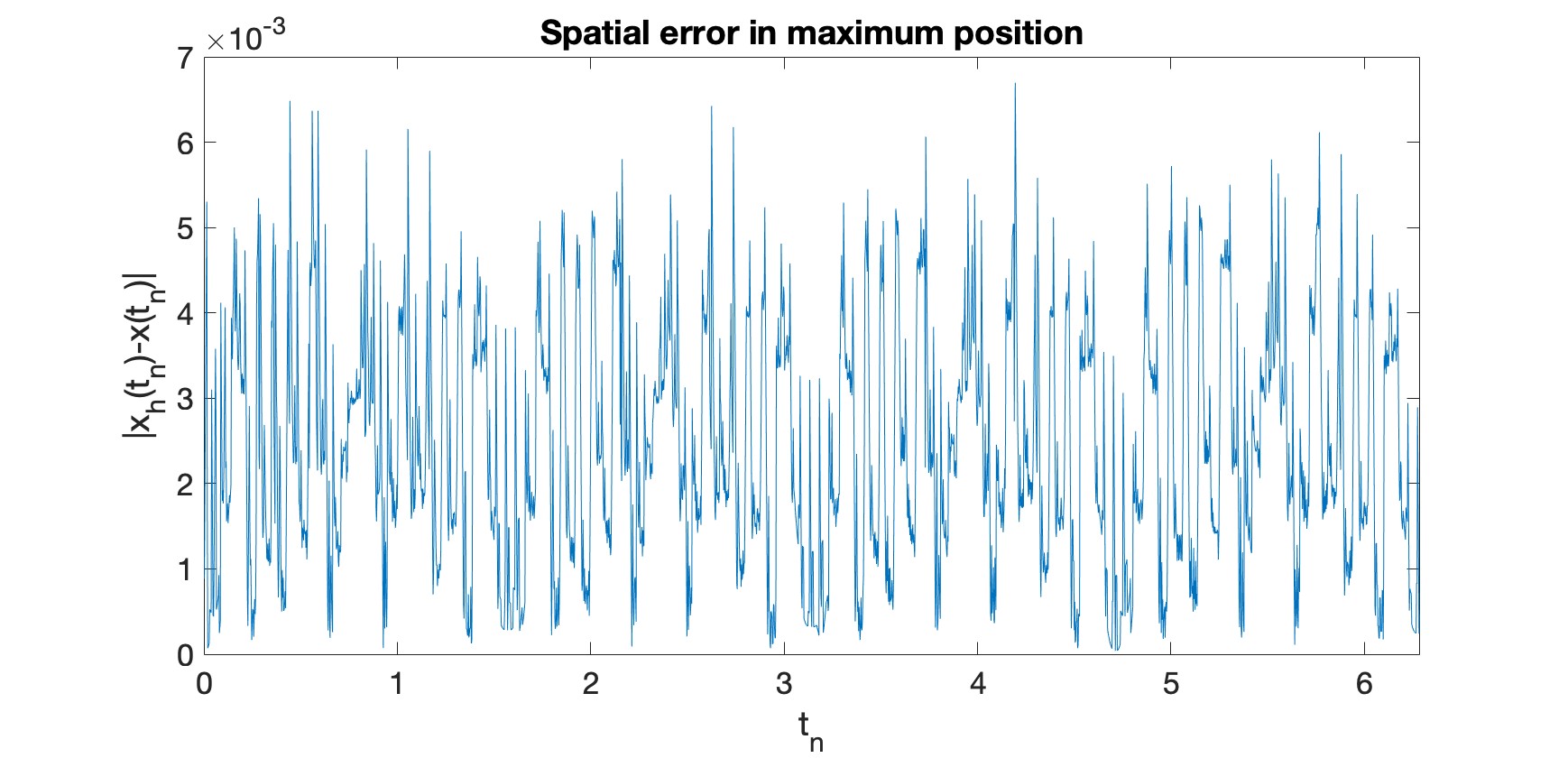}
    \caption{\small\label{F:ex5 maximums comparison} \Change{}{Comparison of the spatial location of the maximum value obtained with the BMZ method and the theoretical position over time. The observed values remain within the range of \(0.00003\) to \(0.00669\), which is notably smaller than the mesh size \(h = 1/80 = 0.0125\). 
    This result highlights the BMZ method's accuracy in tracking the expected position of the maximum with minimal error.}}
\end{figure}

Figure \ref{F:ex5 maximums comparison} presents the spatial distance $|x_h(t_n)-x(t_n)|$ between the location $x_h(t_n)$ of the maximum of \Change{the bilinear part $u_L$ of}{ the discrete solution obtained with BMZ method at each time step and the theoretical maximum position $x(t_n)$, considering the bilinear part $u_L$ and the values of $u_B,u_S$, corresponding to the element and patch bubbles}. The exact position is determined by considering the corresponding pure transport problem, obtained by setting \(\epsilon = 0\). The observed distances range from \(0\) to \(0.007\), which is highly satisfactory given that the mesh size is \( h = \frac{1}{80} = 0.0125 \). This result demonstrates the high accuracy of the BMZ method in capturing the correct transport dynamics and preserving the expected position of the solution's maximum with minimal spatial error. \delete{The observed oscillations stem from the fact that the maximum of $u_L$ is always attained at a vertex of the mesh.}

\section{Conclusions}

We have presented a stabilized algorithm for advection-dominated advection-reaction-diffusion equations, which is based on adding to the usual finite element space residual-free bubbles supported on elements, but also residual-free bubbles supported on patches of two elements; the BMZ method.

The idea of using patch-supported bubbles was first introduced by Cangiani and S\"uli in~\cite{Cangiani2005long,Cangiani2005short}, using a different kind of bubbles, and resorting to Shishkin meshes inside elements for their computation.
Instead, we propose to use a recursive algorithm, using the same bubbles on finer meshes over the elements in order to compute the bubbles and their interactions, which are in turn necessary for the assembling of the outer linear system. This recursion stops after $O(\hspace{0.05 cm} \log (\hspace{0.05 cm} \Change{|\aaa| h \hspace{0.05 cm}/\hspace{0.05 cm}\epsilon}{\Pe_h} ))$ steps and provides very accurate representation of the bubbles, which stabilize the systems and provide strikingly good results.

The computed discrete solutions satisfy the maximum principle for smooth data, without any observable smearing effect.
We have also computed solutions corresponding to discontinuous boundary data, which fall beyond any possible theory, and the results are excellent.

We have computed experimental orders of convergence, which, for the BMZ method, are approximately 2 for the $L^1$ and $L^2$ norms, $1$ for the $H^1$ norm and $1.5$ for the stability norm. Such rates are also observed for the DG method, but this method seems to need some tuning depending on $\epsilon$; the solutions show some spurious oscillations for larger values of $\epsilon$.
On the other hand, smaller rates are observed for the original RFB method, which coincides with SUPG; in some cases the behavior is a bit erratic and the error does not decrease monotonically.

We have also extended our ideas to problems with non-constant coefficients and to the unsteady case. 
Our method demonstrates an outstanding performance across multiple aspects. Only minor negligible oscillations are observed, and when considering the evolution of a pyramid, its edges are very well preserved, and the maximum and minimum values of the solutions align closely with the expected ones.

We have presented our ideas on a square domain, using square elements, in order to computationally test them and assess their performance.
We are convinced that this method can be implemented without too much effort for general meshes; the remaining challenge is just to design convenient data structures and develop a good code.
\Change{Depending on the initial mesh, all the patches that will appear in the recursion will belong to a slightly larger set of equivalence classes, which will not grow after the first recursion step.}{A good starting point for this endeavor could be~\cite{t8code}.} 
It is also worth noticing that the computations corresponding to all the elements of the initial mesh (using the recursion) can be done offline, and moreover, some kind of database can be designed in order to access during the actual computation of the solutions.

Being now convinced that there is a lot of potential in this method, we will continue working on this subject, along the following lines:
\begin{itemize}
    \item More general domains, with triangular as well as quadrilateral meshes.
    \item Three-dimensional domains, with uniform and non-uniform meshes.
    \item Error analysis of the method.
    \item \new{Adaptivity, and variable velocity fields with sharp boundary layers.}
\end{itemize}

\subsection*{Acknowledgments}
The second and third authors were partially supported by Agencia Nacional de Promoci\'on Cient\'ifica y Tecnol\'ogica through grant PICT-2020-SERIE A-03820, and by Universidad Nacional del Litoral through grants CAI+D-2020 50620190100136LI and CAI+D-2024 85520240100018LI.
\new{We would like to thank an anonymous referees for several comments that have helped us to substantially improve this article.}

\bibliographystyle{abbrv}
\bibliography{biblioRFB}
\end{document}